\documentclass[a4paper,11pt]{article}
\usepackage{array}
\usepackage{theorem}
\usepackage{amsmath,amscd,amssymb} 
\usepackage{latexsym}
\usepackage{epsfig}
\usepackage{xypic}
\theorembodyfont{\sl}

\newtheorem{lemma}{Lemma}[section]

\newtheorem{proposition}[lemma]{Proposition}
\newtheorem{theorem}[lemma]{Theorem}
\newtheorem{corollary}[lemma]{Corollary}
\newtheorem{example}[lemma]{Example}
\newtheorem{definition}[lemma]{Definition}
\newtheorem{remark}[lemma]{Remark}
\newtheorem{question}[lemma]{Question}

\newcommand{\CC}{\mathbb C}

\newcommand{\NN}{\mathbb N}

\newcommand{\PP}{\mathbb P}
\newcommand{\QQ}{\mathbb Q}
\newcommand{\RR}{\mathbb R}

\newcommand{\ZZ}{\mathbb Z}

\renewcommand{\cD}{\mathcal D}

\renewcommand{\cL}{\mathcal L}
\newcommand{\cM}{\mathcal M}

\newcommand{\cO}{\mathcal O}

\newcommand{\cX}{\mathcal X}
\newcommand{\cY}{\mathcal Y}

\newcommand{\isoTo}{\stackrel{\sim}{\To}}
\newcommand{\To}{\longrightarrow}
\newcommand{\Mapsto}{\mapstochar\longrightarrow}

\renewcommand{\Tilde}{\widetilde}
\renewcommand{\Hat}{\widehat}
\renewcommand{\Bar}{\overline}

\newcommand{\imic}{\cong}

\newcommand{\sfrac}[2]{{\textstyle{\frac{#1}{#2}}}}

\newcommand{\SL}{\mathop{\mathrm {SL}}\nolimits}

\newcommand{\SU}{\mathop{\mathrm {SU}}\nolimits}
\newcommand{\Orth}{\mathop{\null\mathrm {O}}\nolimits}
\newcommand{\POrth}{\mathop{\null\mathrm {PO}}\nolimits}
\newcommand{\POrthtd}{\mathop{\null\mathrm {P\widetilde{O}}}\nolimits}
\newcommand{\POrthhat}{\mathop{\null\mathrm {P\widehat{O}}}\nolimits}
\newcommand{\PGamma}{\mathop{\null\mathrm {P\Gamma}}\nolimits}

\newcommand{\Aut}{\mathop{\mathrm {Aut}}\nolimits}

\newcommand{\Def}{\mathop{\mathrm {Def}}\nolimits}
\newcommand{\Hilb}{\mathop{\mathrm {Hilb}}\nolimits}
\newcommand{\Mon}{\mathop{\mathrm {Mon}}\nolimits}
\newcommand{\Ref}{\mathop{\mathrm {Ref}}\nolimits}

\newcommand{\vol}{\mathop{\mathrm {vol}}\nolimits}

\newcommand{\mass}{\mathop{\mathrm {mass}}\nolimits}

\newcommand{\disc}{\mathop{\mathrm {disc}}\nolimits}

\newcommand{\id}{\mathop{\mathrm {id}}\nolimits}

\newcommand{\ord}{\mathop{\null\mathrm {ord}}\nolimits}
\newcommand{\gen}{\mathop{\null\mathrm {gen}}\nolimits}

\newcommand{\bdw}{\mathbf w}

\newcommand{\bdN}{\mathbf N}

\newcommand{\gothM}{\mathfrak M}
\newcommand{\sfM}{\mathsf M}
\newcommand{\latt}[1]{{\langle{#1}\rangle}}
\renewcommand{\div}{\mathop{\mathrm {div}}\nolimits}
\newcommand{\Kthree}{\mathop{\mathrm {K3}}\nolimits}
\newcommand{\qedsymbol}{\mbox{$\Box$}}
\newcommand{\qed}{\unskip\nobreak\hfil\penalty50\hskip1em\hbox{}\nobreak
\hfill\qedsymbol\parfillskip=0pt\finalhyphendemerits=0}
\newenvironment{proof}{\begin{ProofwCaption}{Proof}}{\end{ProofwCaption}}
\newenvironment{ProofwCaption}[1]
 {\addvspace\theorempreskipamount \noindent{\it #1.}\rm}
 {\qed \par \addvspace\theorempostskipamount}

\begin{document}

\title{Moduli spaces of irreducible symplectic manifolds}
\author{V.~Gritsenko, K.~Hulek and G.K.~Sankaran}
\maketitle
\begin{abstract}
\noindent
We study the moduli spaces of polarised irreducible symplectic
manifolds. By a comparison with locally symmetric varieties of
orthogonal type of dimension~$20$, we show that the moduli space of
polarised deformation $\Kthree^{[2]}$ manifolds with polarisation of
degree $2d$ and split type is of general type if $d\ge 12$.

\end{abstract}

\section{Introduction}\label{sec:intro}

A simply-connected compact complex K\"ahler manifold is called an
irreducible symplectic manifold if it has an everywhere nondegenerate
$2$-form, unique up to scalar. Irreducible symplectic manifolds are also known as
irreducible hyperk\"ahler manifolds, and for brevity are frequently
referred to simply as symplectic manifolds, omitting the word
``irreducible''. They have been extensively studied by Beauville,
Bogomolov, Debarre, Fujiki, Huybrechts, Markman, Namikawa and O'Grady
among others. Irreducible symplectic manifolds have even complex
dimension: in the surface case they are the $\Kthree$ surfaces.
However, relatively few examples are known. Backgound, and
considerable detail, may be found in Huybrechts' lecture
notes~\cite{Huy3}.

The second cohomology $H^2(X,\ZZ)$ of a symplectic manifold $X$
carries a nondegenerate quadratic form $q_X$ of signature
$(3,b_2(X)-3)$, called the Beauville form, or Beauville-Bogomolov
form. Usually the lattice $L=\big(H^2(X,\ZZ),q_X\big)$ is not
unimodular, nor is it known to be necessarily even, although it is
even in all known examples. A polarisation on $X$ is a choice of ample
line bundle on $X$, or equivalently the cohomology class $h$ of an
ample line bundle. The (Beauville) degree of the polarisation is
defined to be $d=q_X(h)$: it is positive. There is a period map
for symplectic manifolds: the global Torelli theorem, however, is
known to fail in some cases (see~\cite{Deb}, \cite{Nam}). 

Our aim in this paper is to study the moduli of polarised symplectic
manifolds by means of the period map. In Section~\ref{sec:symplectic}
we describe this construction precisely, prove that the moduli spaces
exist and show how they are related to locally symmetric varieties of
orthogonal type: see Theorem~\ref{theo:moduli}. These varieties are
associated with the orthogonal complement $L_{h^\perp}$ of $h$ in $L$. What
lattice $L_{h^\perp}$ is depends in general on the choice of $h$, not just on
the degree as in the case of $\Kthree$ surfaces. The results of Section~\ref{sec:symplectic} hold
for all classes of irreducible symplectic manifolds.

In Section~\ref{sec:monodromy} we specialise to the case of
deformation $\Kthree^{[n]}$ manifolds: that is, symplectic manifolds
deformation equivalent to $\Hilb^n(S)$ for a $\Kthree$ surface~$S$. 
In this case $L=3U\oplus 2E_8(-1)\oplus \latt{-2(n-1)}$. 
Here one has a better understanding of the map from the moduli space
to the locally symmetric variety, thanks to the work of
Markman~\cite{Mar2}. We show in Theorem~\ref{theo:cover} that in this case
one may consider the quotient by the group $\Tilde\Orth(L,h)$ of
automorphisms of $L$ that fix $h$ and act trivially on the
discriminant group $L^\vee/L$. 

To continue further we need to study the orthogonal groups that can
arise. We do this in Section~\ref{sec:groups}, where we mainly study
the lattice $L_{2t}=3U\oplus 2E_8(-1)\oplus \latt{-2t}$. This
leads us to a description of the possible types of polarisation for
deformation $\Kthree^{[n]}$ manifolds. There are two special types,
having only one orbit of polarising vectors. 

For the rest of the paper we are concerned with the case $n=2$ and
with the simplest polarisation, namely the split type (see
Definition \ref{def:split}), where the
lattice $L_{h^\perp}$ is $L_{2,2d}=2U\oplus 2E_8(-1)\oplus \latt{-2} \oplus
\latt{-2d}$. Our main theorem, Theorem~\ref{gentype}, states that
every component of the corresponding moduli space is of general type
as long as $d\ge 12$. O'Grady~\cite{OG4} studied the case $d=1$ and
showed that the moduli space is unirational. There seem to be very 
few other previous results about
dimension~$20$ moduli spaces of orthogonal type. Voisin~\cite{Vo1}
proved that one of them is birational to the moduli space of cubic
fourfolds, and thus unirational, but the type of the polarisation in
that case is not split. Another non-split case (with $d=11$) was
studied recently by Debarre and Voisin~\cite{DV}: again the moduli
space is unirational. In the split case there are only nine possibly
unirational moduli spaces (for $d=9$ and $d=11$ the Kodaira dimension
is non-negative): for polarised $\Kthree$ surfaces there are still
forty-three such possibilities.

The proof of Theorem~\ref{gentype} is similar in style to the
corresponding result for $\Kthree$ surfaces proved in \cite{GHS:K3} 
(see also \cite{Vo2}),
but there are many differences. We use the low-weight cusp form trick,
which guarantees that once the stable orthogonal group
$\Tilde\Orth(L_{2,2d})$ has a cusp form with suitable vanishing of
weight less than the dimension of the moduli space
then the components are of general type. 

We construct the cusp form by means of the quasi pull-back of the
Borcherds form, as in~\cite{GHS:K3}. To do so one requires a vector in
$E_7$ orthogonal to at least $2$ and at most $14$ roots, of
length~$2d$. 

Here there is a significant technical difficulty. The proof that these
vectors exist involves estimating the number of ways of representing
certain integers by various root lattices of odd rank.  In
Theorem~\ref{rtS} we give a new, clear, formulation of Siegel's
formula for this number in the odd rank case. It may be expressed
either in terms of Zagier $L$-functions or in terms of the H.~Cohen
numbers. This analytic estimate shows that the vectors we want exist
for $d\ge 20$, and we can improve this bound slightly by means of a
computer search. 
\smallskip

{\fontsize{9pt}{11pt}\selectfont

  \noindent {\bfseries Acknowledgements}: We should like to thank Eyal
  Markman, Daniel Huybrechts and Kieran O'Grady for informative
  conversations and advice. We also thank the DFG for financial
  support under grant Hu/337-6. The first and the third author are
  grateful to Leibniz Universit\"at Hannover for hospitality. The
  first and the second author would like to thank Max-Planck-Institut
  f\"ur Mathematik in Bonn for support and for providing excellent
  working conditions. We are grateful to the referee for pointing out an error in 
  the first version of the paper.\par}
  
\section{Irreducible symplectic manifolds}\label{sec:symplectic}

In this section we collect the necessary results concerning symplectic
manifolds and their moduli. The main aim is to relate moduli spaces of
polarised symplectic manifolds to quotients of homogeneous domains by
an arithmetic group. 

We begin with the basic definitions and facts about irreducible
symplectic manifolds. 

\begin{definition}\label{def:irrsymplectic}
  A complex manifold $X$ is called an {\it irreducible symplectic
  manifold\/} if the following conditions are fulfilled:
\begin{itemize}
\item[{\rm (i)}] $X$ is a compact K\"ahler manifold;
\item[{\rm (ii)}] $X$ is simply-connected;
\item[{\rm (iii)}] $H^0(X,\Omega^2_X) \cong \CC \sigma$ where $\sigma$
  is an everywhere nondegenerate holomorphic $2$-form. 
\end{itemize}
\end{definition}
Irreducible symplectic manifolds are also known as irreducible
hyperk\"ahler manifolds, and very often simply as symplectic
manifolds. The symplectic surfaces are the $\Kthree$ surfaces. In
higher dimension the known examples are the Hilbert schemes
$\Hilb^n(S)$ of a $\Kthree$ surface $S$ and deformations of them;
generalised Kummer varieties and their deformations; and two examples
of dimensions $6$ and $10$, constructed by O'Grady using moduli spaces
of sheaves on abelian surfaces and $\Kthree$ surfaces respectively
(\cite{OG1}, \cite{OG2}). 

It follows immediately from the definition that $X$ must have even
dimension $2n$ and that its canonical bundle $\omega_X$ is trivial. 
Moreover $h^{2,0}(X)=h^{0,2}(X)=1$ and $h^{1,0}(X)=h^{0,1}(X)=0$. By
a result of Bogomolov \cite{Bog}, the deformation space of $X$ is
unobstructed. This result was generalised to Ricci-flat manifolds by
Tian~\cite{Ti} and Todorov~\cite{Tod}, and algebraic proofs were given
by Kawamata~\cite{Ka} and Ran~\cite{Ran} (see also~\cite{Fuj}). 
Since
\[
T_{[0]} \Def(X) \cong H^1(X,T_X) \cong H^1(X, \Omega^1_X)
\]
the dimension of the deformation space is $b_2(X)-2$. 

The main discrete invariants for symplectic manifolds are the
{\it Beauville form\/} (also known as the {\it Beauville-Bogomolov
 form\/}) and the {\it Fujiki constant\/} or {\it Fujiki invariant}. 
The Beauville form is an indivisible integral symmetric bilinear form
on $H^2(X,\ZZ)$ of signature $(3,b_2(X)-3)$. Its role in the theory of
irreducible symplectic manifolds is similar to the role of the
intersection form for $\Kthree$ surfaces. To define it, let $\sigma
\in H^{2,0}(X)$ be such that $\int_X(\sigma\Bar{\sigma})^n=1$ and
define
\[
q'_X(\alpha)=\frac{n}{2}\int_X\alpha^2(\sigma\Bar{\sigma})^{n-1} +
(1-n)\left(\int_X\alpha\sigma^{n-1}\Bar{\sigma}^{n}\right)
\left(\int_X\Bar{\alpha}\sigma^{n}\Bar{\sigma}^{n-1}\right). 
\] 
After multiplication by a positive constant $\gamma$ the quadratic
form $q_X=\gamma q'_X$ defines an indivisible, integral, symmetric
bilinear form $(\quad ,\quad )_X$ on $H^2(X,\ZZ)$: this is the
Beauville form. Clearly $(\sigma,\sigma)_X=0$ and
$(\sigma,\overline{\sigma})_X >0$. Let $v(\alpha)=\alpha^{2n}$ be
given by the cup product. Then, by a result of Fujiki~\cite[Theorem
4.7]{Fuj}, there is a positive rational number $c$, the {\it Fujiki
invariant\/} such that
\[
v(\alpha)=cq_X(\alpha)^n
\]
for all $\alpha \in H^2(X,\ZZ)$. 

In \cite{OG3} O'Grady introduced the notion of {\it numerical
 equivalence\/} among symplectic manifolds. Two symplectic manifolds
$X$ and $X'$ of dimension $2n$ are said to be numerically equivalent
if there exists an isomorphism $f\colon H^2(X,\ZZ)\isoTo H^2(X',\ZZ)$
of abelian groups with $\int_X \alpha^{2n} = \int_{X'} f(\alpha)^{2n}$
for all $\alpha \in H^2(X,\ZZ)$. The equivalence class of $X$ is
called the {\it numerical type\/} of $X$, denoted by $\bdN$. Clearly,
two symplectic manifolds are numerically equivalent if they have the
same Beauville form and Fujiki invariant. O'Grady \cite[Section
2.1]{OG3} showed that the converse is also true unless
$b_2(X)=b_2(X')=6$ and $n$ is even, in which case the numerical type
determines $c_X$ but {\it a priori\/} one only has $q_{X}= \pm
q_{X'}$. (There are, however, no known examples of irreducible
symplectic manifolds with $b_2=6$.) 

We fix an abstract lattice $L$ which is isomorphic to $H^2(X,\ZZ)$
equipped with the Beauville form $(\quad ,\quad )_X$ (the Beauville
lattice) and consider its associated period domain
\[
  \Omega_L=\{[\bdw] \in \PP(L\otimes \CC) \mid 
  (\bdw,\bdw)=0,\ (\bdw,\Bar\bdw)>0\}
\]
which, since the signature is $(3,b_2(X)-3)$, is connected. A {\it
 marking\/} of a symplectic manifold $X$ is an isomorphism
$\psi\colon H^2(X,\ZZ) \isoTo L$ of lattices. We can associate to each
marked symplectic manifold $(X,\psi)$ its period point $[\psi(\sigma)]
\in \Omega_L$. Now let $f\colon \cX \to U$ be a representative
of $\Def (X)$. This means that $U$ is a polydisc, $X_0:=f^{-1}(0)
\cong X$ and $f$ is a proper submersive map whose Kodaira-Spencer map
\[
T_{f,0}\colon T_{U,0} \To H^1(X,T_X)
\] 
is an isomorphism. We can use the marking $\psi$ to define an
isomorphism $\psi_U: R^2f_*(\ZZ) \isoTo L_U$ (we shall tacitly shrink
$U$ wherever necessary) and thus a period map
\begin{align*}
\varphi_U\colon U &\To \Omega_L\\
t &\Mapsto [\psi_t (\sigma_{X_t})]. 
\end{align*}
The local Torelli theorem for symplectic manifolds, proved 
by Beauville \cite{Be}, says that $\phi_U$ is a local isomorphism (in
the complex topology). 

The surjectivity of the period map was proved by Huybrechts in
\cite[Theorem 8.1]{Huy1}. To formulate his result we consider a fixed
lattice $L$ which appears as the Beauville lattice of some symplectic
manifold. Let $\gothM_L$ be the corresponding moduli space of marked
symplectic manifolds, {i.e.} as a set $\gothM_L=\{(X,\psi\colon H^2(X,
\ZZ) \isoTo L \} / \approx $ where the equivalence relation $\approx$ is
induced by $\pm f^*$ with $f: X \to X'$ a biholomorphic map. The space
$\gothM_L$ admits a natural smooth complex structure which, however, is
not Hausdorff. The period map $\varphi\colon \gothM_L \to \Omega_L $
is a holomorphic map and Huybrechts has shown that every connected
component of $\gothM_L$ maps surjectively onto $\Omega_L$. 
For a discussion of moduli of marked symplectic manifolds see the
paper~\cite{Huy5} by Huybrechts. 

The situation improves considerably when one considers moduli of
polarised symplectic manifolds. A {\it polarisation\/} on a symplectic
manifold $X$ is the choice of an ample line bundle $\cL$ on $X$. 
Since the irregularity of $X$ is $0$ this is the same as the choice of
a class $h \in H^2(X,\ZZ)$ representing an ample line bundle on $X$. 
Clearly $q_X(h) >0$. Conversely, Huybrechts has shown (\cite[Theorem
3.11]{Huy1}: see also \cite[Theorem 2]{Huy2}) that a symplectic
manifold $X$ is projective if and only if there exists a class $h \in
H^2(X,\ZZ) \cap H^{1,1}(X)$ with $q_X(h) > 0$. It should be noted,
however, that neither line bundle associated to $\pm h$ need be ample. 
There is, however, a small deformation of the pair $(X,h)$ with this
property. 

We now fix an abstract lattice $L$ of rank $b_2 = b_2(X)$ such that
$H^2(X,\ZZ) \cong L$ and let $h \in L$ be a primitive element
with $h^2 >0$. Then the lattice
\[
L_{h^\perp} = h^{\perp}_L < L
\]
has signature $(2,b_2-3)$. It defines a homogeneous domain, which in
this case has two connected components
\begin{equation}\label{cD}
\Omega_{L_{h^\perp}} = \cD(L_{h^\perp}) \cup \cD'(L_{h^\perp}). 
\end{equation}
If $(X,h)$ is a pair with $h \in H^2(X,\ZZ) \cap H^{1,1}(X)$ and
$\psi\colon H^2(X,\ZZ) \to L$ is a marking then the period point
$[\psi(\sigma)] \in \Omega_{L_{h^\perp}}$. Hence for every deformation $\cX
\to U$ of the pair $(X,h)$ the period map defines a holomorphic map
$\varphi_U\colon U \to \Omega_{L_{h^\perp}}$. 

In this paper we are interested in the moduli spaces of polarised
symplectic manifolds. We shall fix the dimension $2n$ and the
numerical type $\bdN$ of the symplectic manifolds that we consider. We
have already remarked that this determines the Beauville lattice and
Fujiki invariant unless $b_2=6$ and $n$ is even, in which case the
quadratic form is only determined up to sign. We shall consider
polarised symplectic manifolds $(X,h)$ of fixed numerical type and
given value $q_X(h)=d >0$. The degree of the associated line bundle
$\cL$ is $\deg(\cL)=h^{2n}=cd^n$ where $c$ is the Fujiki invariant. 
Instead of working with the (geometric) degree of a polarisation we
prefer to work with the number $d$, which we will call the {\it
Beauville degree\/} of the polarisation. We first note the following
variant of a result of Huybrechts~\cite[Theorem 4.3]{Huy4}, which is
itself an application of  the finiteness theorem of Koll\'ar and Matsusaka
\cite[Theorem 3]{KM}. 

\begin{proposition}\label{prop:finiteness}
  For fixed numerical type there are only finitely many deformation
  types of polarised symplectic manifolds $(X,h)$ of dimension $2n$
  and given Beauville degree $d=q_X (h) > 0$. 
\end{proposition}
\begin{proof}
  Since the numerical type determines the Fujiki invariant $c$ our
  choices also fix the degree $h^{2n} = cq_X (h)^n > 0$. The result
  follows immediately from~\cite[Corollary 26.17]{Huy4}. 
\end{proof}

Now we define the moduli spaces we are interested in. We first fix a
possible Hilbert polynomial, say $P(m)$.  Note that this is more than
fixing the degree of the polarisation.  By Matsusaka's big theorem we
can find a constant $m_0$ such that for all polarised manifolds
$(X,\cL)$ with Hilbert polynomial $P(m)$ the line bundles
$\cL^{\otimes m}$ are very ample for $m \geq m_0$ and have no higher
cohomology.  Then we have embeddings $\varphi_{|\cL^{\otimes
 m_0}|}\colon X \to \PP^{N-1}$ where $N=h^0(X,\cL^{\otimes
 m_0})=P(m_0)$.  Such an embedding depends on the choice of a basis
of $H^0(X,\cL^{\otimes m_0})$.  Let $H$ be an irreducible component of
the Hilbert scheme $\Hilb_P(\PP^{N-1})$ that contains at least one
point $\eta\in H$ corresponding to a symplectic manifold $X_\eta$. We
denote by $H_{\text{sm}}$ the open part of $H$ parametrising smooth varieties. 
The following lemma is well known. We include a proof since we are not aware of a suitable reference.

\begin{lemma} \label{lem:hilbert}
$H_{\text{sm}}$ has the following properties:
\begin{itemize}
\item[\rm (i)] Every point in $H_{\text{sm}}$ parametrises a symplectic manifold;
\item[\rm (ii)] $H_{\text{sm}}$ is smooth. 
\end{itemize}
\end{lemma}
\begin{proof}
Claim (i)~follows from
 Beauville's classification theorem \cite[\S 5, Th\'eor\`eme 2]{Be}:
 the universal family $\cX_{\text{sm}}$ over $H_{\text{sm}}$ is a flat family of
  projective (and hence compact K\"ahler) manifolds which are simply-connected with trivial canonical bundle.
Moreover, since the second
  Betti number is constant $h^{2,0}(X)=1$ for
  every for every fibre by semi-continuity. 

This can be proved along the lines of \cite[Theorem 1.3]{Sz}.
The long exact sequence of the normal bundle sequence yields
\[
\cdots \To H^0(X,N_{X/\PP^{N-1}}) \stackrel{\alpha} \To H^1(X,T_X) \To
H^1(X,T_{\PP^{N-1}}|_X) \To \cdots 
\]
The image of $\alpha$ is contained in the hyperplane $V_h=h^{\perp}
\subset H^1(X,T_X) \cong H^1(X,\Omega^1_X)$, which corresponds to
deformations of the pair $(X,h)$ where $h=c_1(\cL)=\cO_X(1)$. 
Since $H^1(X,T_{\PP^{N-1}}|_X)$ is $1$-dimensional (which follows from
the restriction of the Euler sequence to $X$) the image
of $\alpha$ is equal to $V_h$ and hence the Hilbert scheme is
unobstructed and thus smooth. 
\end{proof}

\begin{definition}\label{def:orthogonals}
Let $L$ be a lattice.
The {\it stable orthogonal group\/}
$\Tilde\Orth(L)$ is defined by
\begin{equation}\label{defOtildeL}
\Tilde{\Orth}(L) = \{ g \in \Orth(L) \mid g(l^{\vee}) \equiv
l^{\vee} \mod L \text{ for all }
l^{\vee} \in L^{\vee}\}.
\end{equation}
We shall also need
\begin{equation}\label{defOhatL}
\Hat{\Orth}(L) = \{ g \in \Orth(L) \mid g(l^{\vee}) \equiv
\pm l^{\vee} \mod L \text{ for all }
l^{\vee} \in L^{\vee}\}. 
\end{equation}
For a primitive element $h \in L$ with $h^2=d >0$, we define the groups
\begin{equation}\label{defOLh}
\Orth(L,h)=\{ g \in \Orth(L) \mid g(h)=h \}
\end{equation}
and \begin{equation}\label{defOtildeLh}
\Tilde{\Orth}(L,h) = \{ g \in \Tilde{\Orth}(L) \mid g(h)=h \}.
\end{equation}
For any subgroup $\Gamma \subset \Orth(L)$ we
define the projective group
$\PGamma = \Gamma / (\pm 1)$. If $L$ is indefinite we define $\Gamma^+$
to be the subgroup of $\Gamma$ of elements with real spin norm~$1$.
\end{definition}
There are two ways to choose the definition of the spin norm. We have
chosen it in such a way that any $-2$-reflection has spin norm~$1$. 

We can consider $\Orth(L,h)$ and $\Tilde{\Orth}(L,h)$ as subgroups of
$\Orth(L_{h^\perp})$, where $L_{h^\perp}$ is, as usual, the lattice perpendicular to
$h$ in $L$.  We shall discuss the relationships among these three
groups in Section~\ref{sec:groups}. 

Note that in our case the lattice $L$ has signature $(3,b_2-3)$ and
hence $L_{h^\perp}$ has signature $(2,b_2-3)$.  Thus the lattice $L_{h^\perp}$
determines a homogeneous domain $\Omega_h=\Omega_{L_{h^\perp}}$ of type~IV on
which the three groups act, and if $\Gamma<\Orth(L_{h^\perp})$ then $\Gamma^+$
is the subgroup of $\Gamma$ that preserves the component~$\cD$. 

The following theorem is crucial for the rest of the paper.  
Viehweg's results on moduli spaces \cite{Vi} together with
Proposition \ref{prop:finiteness} give us the existence of moduli spaces
of polarised irreducible symplectic manifolds of fixed numerical type.
More importantly, this result allows us to relate these moduli spaces to modular
varieties, in this case to quotients of a homogeneous domain of type IV by a suitable
arithmetic group. (Both the domain and the group depend on the moduli problem
in question.)

\begin{theorem}\label{theo:moduli}
  There exists a quasi-projective coarse moduli space
  $\sfM_{2n,\bdN,d}$ parametrising primitively polarised symplectic
  manifolds of dimension $2n$, numerical type $\bdN$ and Beauville
  degree~$d$. We choose any one of the irreducible components of
  $\sfM_{2n,N,d}$ and denote it by $\cM_d$. Such a choice determines
  a primitive vector $h\in L$ (or possibly $h\in L(-1)$ if $b_2=6$ and
  $n$ is even) with $q(h)=d$ such that there is a map
\[
\varphi\colon  \cM_d \To \left(\Orth(L,h) \setminus
  \Omega_h \right)^0.  
\]
Here $(\Orth(L,h) \setminus \Omega_h)^0$ is a connected
component of $\Orth(L,h) \setminus \Omega_h$. 
The map $\varphi$ is a morphism of quasi-projective varieties
which is dominant and finite onto its image. 
\end{theorem}
\begin{proof}
  Proposition \ref{prop:finiteness} shows that there are
  only finitely many possible Hilbert polynomials for a given choice
  of the discrete data $2n$, $\bdN$ and $d$. It then follows from Viehweg's work
  \cite[Theorem 1.13]{Vi} and the discussion there that the moduli spaces in question exist: 
  indeed, every component $\cM_d$ of $\sfM_{2n,\bdN,d}$ is a
  quotient of the form $\SL(N,\CC) \backslash H_{\text{sm}}$ for some
  component $H$ of a suitable Hilbert scheme (see the discussion of
  Lemma~\ref{lem:hilbert}). 

  We now want to relate the components $\cM_d$ to quotients of the
  form $\Orth(L,h) \backslash \Omega_h$.  For this we want to
  construct a map $\tilde\varphi\colon H_{\text{sm}} \to \Orth(L,h) \backslash
  \Omega_h$ and then argue that it factors through the quotient by
  $\SL(N,\CC)$. We first observe that every component $H$ determines
  an $\Orth(L)$-orbit of primitive vectors $h \in L$ with $q(h)=d$. 
  Indeed, choosing a local marking $\psi_t$ near a given point in
  $H_{\text{sm}}$ we obtain a vector $h_t=\psi_t\big(c_1(\cO_{X_t}(1))\big)$ with
  $q(h_t)=d$, and any two local markings differ by an element of
  $\Orth(L)$.  Since $H_{\text{sm}}$ is connected and the number of
  $\Orth(L)$-orbits is finite this associates to each $H_{\text{sm}}$ a unique
  $\Orth(L)$-orbit. 
  Let $h$ be a representative of the orbit defined by $H_{\text{sm}}$. We shall
  be interested only in {\it $h$-markings\/}, that is, markings
  $\psi$ with $\psi\big(c_1(\cO_X(1))\big)=h$. They exist locally on
  all of $H_{\text{sm}}$, and an $h$-marking on an open set $U\subset H_{\text{sm}}$
  defines, via the period map, a holomorphic map $\varphi_U\colon U
  \to \Omega_h$. Two $h$-markings differ by an element $\Orth(L,h)$,
  so we obtain a holomorphic map $\tilde\varphi\colon H_{\text{sm}} \to
  \Orth(L,h) \backslash \Omega_h$. 

  Assume that $M \in \SL(N,\CC)$ maps $(X,\cO_X(1))$ to $(X',\cO_{X'}(1))$.
Given $h$-markings $\psi\colon H^2(X,\ZZ) \to L$ and $\psi' \colon H^2(X',\ZZ) \to L$
there exists an element $g \in \Orth(L,h)$ with
$\psi \circ M^* = g \circ \psi'$. This shows that the map
$\tilde\varphi$ factors through the quotient by $\SL(N,\CC)$, giving
the required map
\[
\varphi\colon  \cM_d \To \left(\Orth(L,h) \backslash
  \Omega_h \right)^0. 
\] 

Next we want to show that the map $\varphi$ is a morphism of
quasi-algebraic varieties. This can be deduced from a theorem of Borel
\cite{Bl} which says the following: if $Y$ is a quasi-projective
variety and $f\colon Y \to \Gamma \backslash \Omega$ a holomorphic map
to an arithmetic quotient of a homogeneous domain, where $\Gamma$ is
torsion free, then $f$ is a morphism of quasi-projective varieties. 
Here $\Gamma \backslash \Omega$ carries the natural structure as a
quasi-projective variety, which comes from the Baily-Borel
compactification. 
Although one cannot apply this theorem immediately, as
$\Orth(L,h)$ will in general not be torsion free, the difficulty
can be avoided by using level structures. The argument given 
by Hassett in \cite[Proposition 2.2.2]{Ha}
in the special case of cubic fourfolds
carries over without difficulties.

Finally we want to prove that $\varphi$ is dominant and has finite
fibres. Since $\varphi$ is a morphism of quasi-projective varieties it is enough
to show that
$\varphi$ has no positive-dimensional fibres and that the image has the same dimension as
the period domain. We proceed by following Szendr\H{o}i's arguments in \cite{Sz}.
As in \cite[Lemma 2.7]{Sz} one can construct a finite \'etale
covering $H'_{\text{sm}} \to H_{\text{sm}}$ with the property that the action of $\SL(N,\CC)$
lifts to a free action on $H'_{\text{sm}}$ as well as to an action on the pullback
$\cY_{\text{sm}} \to H'_{\text{sm}}$ of the universal family $\cX_{\text{sm}} \to H_{\text{sm}}$.
Dividing by the action of $\SL(N,\CC)$ one obtains a quotient family
over a base $Z_{\text{sm}}$ which is smooth and finite over the moduli space.
By the infinitesimal Torelli theorem the
period map on $Z_{\text{sm}}$ is open with discrete fibres near every point of $Z_{\text{sm}}$. 
Since the group $\Orth(L,h)$ acts properly discontinuously on the
domain $\Omega_h$, the induced morphism $Z_{\text{sm}} \to \Orth(L,h)
\backslash \Omega_h$ is dominant and has no positive-dimensional
fibres and hence the same also holds for $\varphi$.
\end{proof}

\begin{remark}\label{rk:nonsurjective}
  The map $\varphi\colon \cM_d \to \left(\Orth(L,h)
  \backslash \Omega_h \right)^0$ will in general not be
  surjective as there are period points in $\Omega_h$ which
  parametrise pairs $(X,h)$ where $h$ is not ample. 
\end{remark}
This phenomenon already occurs for $\Kthree$ surfaces. Unlike in the
$\Kthree$ case it is, however, not clear which open part of the period
domain belongs to ample divisors. There are some results about this in
special cases, due to Hassett and Tschinkel~\cite{HT}. 

We shall use Theorem~\ref{theo:moduli} in Section~\ref{sec:modular} to
prove general type results for some moduli spaces of symplectic
manifolds by proving that the quotients $\Orth(L,h) \backslash
\Omega_h$ are of general type. 

\section{Deformation $\Kthree^{[n]}$ manifolds and
  monodromy}\label{sec:monodromy}

For the remainder of the paper we concentrate on a special case. 
\begin{definition}
  A {\it deformation $\Kthree^{[n]}$ manifold\/} is a symplectic
  manifold that is deformation equivalent to $\Hilb^n(S)$ for some
  $\Kthree$ surface~$S$. 
\end{definition}
(Compare the definition of numerical $\Kthree^{[2]}$ in~\cite{OG3}.) If
$X$ is a deformation $\Kthree^{[n]}$ manifold then $H^2(X,\ZZ)\imic
L_{2n-2}$ (as a lattice with the Beauville form), where for any $t\in
\NN$ we put
\begin{equation}\label{defineL2t}
L_{2t} = 3U \oplus 2E_8(-1) \oplus \latt{-2t}. 
\end{equation}
For deformation $\Kthree^{[n]}$ manifolds, the numerical type is
determined completely by the dimension $2n$, and the (Beauville)
degree of a polarisation is always even. The Fujiki invariant is
$\frac{(2n)!}{2^n n!}$. 

We study deformation $\Kthree^{[n]}$ manifolds by using monodromy
operators, whose theory was developed by Markman \cite{Mar1},
\cite{Mar2}. We consider a flat family $\pi\colon \cX \to
B$ of compact complex manifolds with fibre $X$ over the point $b \in
B$.  Associated to such a family we obtain a monodromy representation
\[
\pi_1(B,b) \To \Aut (H^*(X,\ZZ)). 
\]  
We define the {\it group of monodromy operators\/} to be the subgroup
$\Mon(X)$ of $\Aut(H^*(X,\ZZ))$ generated by the image of all
monodromy representations. If we restrict to the second cohomology, we
obtain a representation $\pi_1(B,b) \to \Aut (H^2(X,\ZZ))$ and
correspondingly a subgroup $\Mon^2(X) \subset \Aut(H^2(X,\ZZ))$.  If
$X$ is a symplectic manifold then monodromy transformations preserve
the Beauville form and we obtain a subgroup $\Mon^2(X) \subset
\Orth(H^2(X,\ZZ))$. 

Let $\Ref(X)$ be the subgroup of $\Orth(H^2(X,\ZZ))$ generated by
$-2$-reflections and by the negatives of $+2$-reflections. By the
choice of spin norm made in Definition~\ref{def:orthogonals}, this is
a subgroup of $\Orth^+(H^2(X,\ZZ))$. 

\begin{theorem}\label{monodromygroup} {\bf (Markman \cite[Theorem
    1.2]{Mar2}.)} 
If $X$ is a deformation $\Kthree^{[n]}$ manifold then
\[
\Mon^2(X) = \Ref(X). 
\]
\end{theorem}
Using a marking $\psi\colon H^2(X,\ZZ) \isoTo L$ we can think of
$\Mon^2(X)$ as a subgroup of $\Orth^+(L_{2n-2})$. Since $\Mon^2(X)$ is a
normal subgroup we obtain a well-defined subgroup
$\Mon^2(L_{2n-2})=\Ref(L_{2n-2})=\Ref(X) \subset \Orth^+(L_{2n-2})$. 

It follows from a result of Kneser \cite[Satz 4]{Kn}
or from \cite[Lemma 4.10]{Mar1} that the
groups satisfy
\begin{equation}\label{reflectiongroup}
\Ref(L_{2n-2}) = \Hat\Orth^+(L_{2n-2})
\end{equation}
(see Definition~\ref{def:orthogonals}).
Note that the assumptions of Kneser's theorem are fulfilled since
$L_{2n-2}$ contains three copies of $U$. 

Unlike in the case of $\Kthree$ surfaces, for fixed degree $2d$ there
is not a unique $\Orth^+(L_{2n-2})$-orbit of primitive vectors $h$ with
$h^2=2d$. We shall address this question in Section~\ref{sec:groups}. 
Hence the moduli space of deformation $\Kthree^{[n]}$ mani\-folds with a
primitive polarisation of degree $2d$ will in general have more than
one component. 

\begin{theorem}\label{theo:cover} 
  Let $\cM^{[n]}_{2d}$ be an irreducible component of the moduli space
  of deformation $\Kthree^{[n]}$ manifolds with a primitive
  polarisation of degree $2d$. Then the map $\varphi$ from Theorem
  \ref{theo:moduli}, above, factors through the finite cover
  $\Tilde\Orth^+(L_{2n-2},h)\backslash\cD_h\to\Orth^+(L_{2n-2},h)\backslash\cD_h$:
  that is, there is a commutative diagram
\[
\xymatrix{
{\cM^{[n]}_{2d}} \ar[r]^(.3){\tilde\varphi} \ar[dr]^{\varphi}
& {\Tilde\Orth^+(L_{2n-2},h) \backslash \cD_h}  
\ar[d]\\
& {\Orth^+(L_{2n-2},h) \backslash \cD_h} 
}
\]
\end{theorem}
\begin{proof}
  Recall from the proof of Theorem \ref{theo:moduli} that
  $\cM^{[n]}_{2d}=\SL(N,\CC) \backslash H_{\text{sm}}$ for some suitable open
  part of a component $H$ of the Hilbert scheme. We choose a base
  point in $H_{\text{sm}}$ and denote the corresponding symplectic variety by
  $X_0$. Choose an $h$-marking $\psi_0\colon
  \big(H^2(X_0,\ZZ),c_1(\cO_{X_0(1)})\big) \to (L_{2n-2},h)$. Now let
  $Y$ be a variety corresponding to another point in $H$ and choose a
  path $\sigma_Y$ from $X_0$ to $Y$. Transporting the marking $\psi_0$
  along this path we obtain an $h$-marking $\psi_{\sigma_Y}\colon
  \big(H^2(Y,\ZZ),c_1(\cO_Y(1))\big) \to (L_{2n-2},h)$.  Clearly this
  marking will depend on the path $\sigma_Y$. Let $\tau_Y$ be another
  path from $X_0$ to $Y$ and $\psi_{\tau_Y}$ the corresponding
  marking.  Then $\tau_Y \circ \sigma_Y^{-1}$ is a closed path based
  at $Y$ and induces an automorphism $f^*= \psi_{\sigma_Y}^{-1} \circ
  \psi_{\tau} \in \Mon^2(Y)$.  Let $f'= \psi_{\sigma_Y} \circ f^*
  \circ \psi_{\sigma_Y}^{-1} \in \Mon^2(L_{2n-2})$.  Then $f' \circ
  \psi_{\sigma_Y} = \psi_{\sigma_Y} \circ f^* = \psi_{\tau_Y}$. This
  shows that we have a morphism
\[
\varphi'\colon H_{\text{sm}} \To (\Mon^2(L_{2n-2}) \cap \Orth^+(L_{2n-2},h))
  \backslash \cD_h = \Hat\Orth^+(L_{2n-2},h) \backslash
  \cD_h
\]
where the last equality follows from Theorem~\ref{monodromygroup} and
equation~(\ref{reflectiongroup}). 

We next claim that $\varphi'$ factors through $\cM^{[n]}_{2d}$.  For
this let $g \in \SL(N,\CC)$ be an element which maps $Y$ to $Z$.  Let
$\sigma_Y$ and $\sigma_Z$ be paths from $X_0$ to $Y$ and $Z$
respectively, with corresponding markings $\psi_{{\sigma_Y}}$ and
$\psi_{{\sigma_Z}}$. We now consider the path $\sigma_Z \circ
\sigma_Y^{-1}$ from $Y$ to $Z$. Using the element $g$ to identify $Y$
and $Z$ makes this a closed path. We can now argue as above and
conclude that $\psi_{{\sigma_Y}} \circ g^* \circ
{\psi_{{\sigma_Z}}}^{-1} \in \Mon^2(L_{2n-2}) \cap \Orth^+(L_{2n-2},h)
=\Hat\Orth^+(L_{2n-2},h)$.  (Strictly speaking we need a
complex family to argue that this element is in $\Mon^2(L_{2n-2})$, but this
can easily be achieved by a complex thickening of the closed path.) 

We can put  $\Tilde\Orth^+(L_{2n-2}, h)$ in the formulation of the theorem
because $\POrthtd^+(L_{2n-2}, h)\imic \POrthhat^+(L_{2n-2}, h)$,
and the groups act on the symmetric space through their projectivisations. 
\end{proof}

\begin{remark}\label{rk:lifts}
  The lifting of the map $\varphi$ to $\tilde\varphi$ is not unique.
  Two markings $\psi_0$ and $\psi_1$ define the same lifting if and
  only if $\psi_0 \circ {\psi_1}^{-1}$ is trivial in
  $\POrth^+(L_{2n-2},h)/\POrthtd^+(L_{2n-2},h)$, so the different
  liftings are classified by the quotient
  $\POrth^+(L_{2n-2},h)/\POrthtd^+(L_{2n-2},h)$. We shall compute the
  index of $\Tilde\Orth^+(L_{2n-2},h)$ in $\Orth^+(L_{2n-2},h)$ below
  (Proposition~\ref{ogroup}), in almost all cases.
\end{remark}
Theorem~\ref{theo:cover} should also be compared to Markman's
consideration of the non-polarised case in \cite[Section 4.2]{Mar2}. 
\begin{remark}\label{rk:no_Torelli}
  As in \cite{Mar2} we can conclude from Theorem~\ref{theo:cover} that
  the global Torelli theorem for polarised deformation $\Kthree^{[n]}$
  manifolds fails whenever $[\POrth^+(L_{2n-2},h) :
  \POrthtd^+(L_{2n-2},h)] > 1$. This can occur: see
  Proposition~\ref{ogroup}, below. 
\end{remark}

With Remark~\ref{rk:no_Torelli} in mind we pose the
following question.  
\begin{question}\label{qu:irreducible}
Is it true that for every $\Orth^+(L_{2n-2})$-orbit
of some primitive vector $h$ with $h^2=2d>0$ the part of the moduli
space $\sfM_{2n,\bdN,2d}$ corresponding to polarisations in the orbit
of $h$ is irreducible and that the map $\tilde\varphi$ has degree~$1$? 
\end{question}
A positive answer to both parts of Question~\ref{qu:irreducible} could
be viewed as the correct version of the global Torelli theorem for
deformation $\Kthree^{[n]}$ manifolds. 

\begin{remark}\label{rk:moninstead}
  For every class of symplectic manifolds, Theorem \ref{theo:moduli}
  remains true if we consider the monodromy group instead of the
  orthogonal group. 
\end{remark}

\section{Orthogonal groups}\label{sec:groups}

Let $L$ be an even lattice. By lattice (or sublattice) we always mean a
non-degenerate lattice (or sublattice). 
We denote the discriminant group of $L$
by $D(L)=L^{\vee}/L$. This carries an induced quadratic form with values in $\QQ / 2\ZZ$.
If $g\in \Orth(L)$ we denote
by $\bar g$ its image in $\Orth(D(L))$. 

Let $S$ be a primitive sublattice of $L$: we are mainly interested 
in the case $S=L_{h^\perp}$ for some $h\in L$ with $h^2\neq 0$, but we want 
to consider this more general situation. Analogously to 
Definition~\ref{def:orthogonals} we define the groups
\[
\Orth(L,S)=\{ g\in \Orth(L)\mid g|_S\in \Tilde\Orth(S)\}
\quad\text{and}\quad \Tilde\Orth(L,S)=\Orth(L,S)\cap \Tilde\Orth(L). 
\]
Note that $\Orth(L,\ZZ h)=\Orth(L, h)$ if $h^2\ne \pm 2$. 

Let $S^{\perp}$ be the orthogonal complement of $S$ in $L$.  We have
\[
S^{\perp}\oplus S< L < L^\vee < (S^\perp)^\vee\oplus S^\vee. 
\]
The overlattice $L$ is defined by the finite  subgroup
\[
H=L/(S^{\perp}\oplus S)<(S^\perp)^\vee/S^\perp \oplus S^\vee/S=
D(S^\perp)\oplus D(S)
\]
which is an isotropic subgroup of $D(S^\perp)\oplus D(S)$.  Following
\cite{Nik} we consider the projections
\begin{alignat*}{2}
p_S\colon H \to D(S),\qquad & p_{S^\perp}\colon H\to D(S^\perp). 
\end{alignat*}
Using the definitions and the fact that the lattices $S$ and
$S^\perp$ are primitive in $L$ one can show (see
\cite[Prop. 1.5.1]{Nik}) that these projections are injective and
moreover that if $d_S\in p_S(H)$ then there is a unique $d_{S^\perp}\in
p_{S^\perp}(H)$ such that $d_S+d_{S^\perp}\in H$. 
\begin{lemma}\label{lem:extension}
$\alpha\in \Orth(S^\perp)$ can be extended to $\Orth(L)$
if and only if
\[
\bar \alpha(p_{S^\perp}(H))=p_{S^\perp}(H)
\]
and there exists $\beta\in \Orth(S)$ such that
$p_S^{-1}\circ\bar{\beta}\circ p_S=
p_{S^\perp}^{-1}\circ\bar{\alpha}\circ p_{S^\perp}$. 
\end{lemma}
This is a reformulation of \cite[Corollary 1.5.2]{Nik}. The following
is a particular case. 

\begin{lemma}\label{ols} 
  Let $S$ be a primitive sublattice of an even lattice $L$. 
\begin{itemize}
\item[{(i)}] $g\in \Orth(L,S)$ if and only if $g(S)=S$, $\bar
  g|_{D(S)}=\id$ and $\bar g|_{p_{S^\perp}(H)}=\id$. 

\item[{(ii)}] $\alpha\in \Orth(S^\perp)$ can be extended to
  $\Orth(L,S)$ if and only if $\bar\alpha|_{p_{S^\perp}(H)}=\id$. 

\item[{(iii)}] If $p_{S^\perp}(H)=D(S^\perp)$ then
  $\Orth(L,S)|_{S^\perp}\cong \Tilde\Orth(S^\perp)$. 

\item[{(iv)}] Assume that the projection $\Orth(S^\perp)\to
  \Orth(D(S^\perp))$ is surjective. Then
\[
\Orth(L,S)|_{S^\perp}/\Tilde\Orth(S^\perp)\cong
\{\bar\gamma \in \Orth(D(S^{\perp}))\mid \bar\gamma|_{p_{S^\perp}(H)}=\id\}. 
\]
\end{itemize}
\end{lemma} 

\begin{remark}\label{rk:duals}
  Let $g\in \Orth(L,S)$.  Then $\bar g|_{p_{S^\perp}(H)}=\id$ is
  equivalent to $\bar g|_H=\id$ or to $\bar g|_{H^\vee}=\id$, where
  $H^\vee=((S^\perp)^\vee\oplus S^\vee)/L^\vee$. 
  The condition $\bar g|_{H^\vee}=\id$ is equivalent to the
  following: for any $v\in (S^\perp)^\vee$ we have $g(v)-v\in
  (S^\perp)^\vee\cap L^\vee$.  But $(S^\perp)^\vee\cap L^\vee$ might
  be larger than $S^\perp$.  This shows in terms of the dual lattices
  that $\bar g|_{p_{S^\perp}(H)}=\id$ is weaker than $g|_{S^\perp}\in
  \Tilde \Orth(S^\perp)$. 
\end{remark}

\begin{corollary}\label{cor:perps}
  If $|H|=|\det S^\perp|$ then $\Orth(L,S)|_{S^\perp}\cong
  \Tilde\Orth(S^\perp)$. 
\end{corollary}
\begin{proof} 
  This follows from the injectivity of $p_{S^\perp}$ on $H$ and from
  Lemma~\ref{ols}(iii).  For example, the condition of
  the corollary is true if $L$ is an even unimodular lattice and $S$
  is any primitive sublattice of $L$. 
\end{proof}

If $l\in L$ its {\it divisor\/} $\div(l)$ is the positive generator of the
ideal $(l,L)\subset \ZZ$. Therefore $l^*=l/\div(l)$ is a primitive
element of the dual lattice $L^\vee$ and $\div(l)$ is a divisor of
$\det (L)$.  We recall the the following classical criterion of
Eichler (see~\cite[\S 10]{E}). 
\begin{lemma}\label{Eichler}
Let $L$ be a lattice containing two orthogonal isotropic planes. Then
the $\Tilde\Orth(L)$-orbit of a primitive vector $l\in L$ is
determined by two invariants: by its length $l^2=(l,l)$ and its image
$l^*+L$ in the discriminant group $D(L)$
\end{lemma}
We consider the special lattice $L_{2t}=3U\oplus 2E_8(-1)\oplus
\latt{-2t}$ defined in Equation~\eqref{defineL2t} above. We shall need
this for the application in Section~\ref{sec:modular}.  It has
signature $(3,20)$. We denote a generator of the $1$-dimensional
sublattice $\latt{-2t}$ by $l_t$, so $l_t^2=-2t$, and we denote by
$h_d$ a primitive vector of length $2d$. Note that $\div(h_d)$ is a
common divisor of $2d$ and $2t=-\det(L_{2t})$.

We now investigate the groups $\Orth(L_{2t}, h_d)$ and
$\Tilde\Orth(L_{2t},h_d)$.  In Proposition~\ref{hdorbits} we count the
$\Tilde\Orth(L_{2t})$-orbits of the polarisation vectors $h_d$. This
number is not always~$1$. Thus the situation here is different from
the case of $\Kthree$ surfaces, because there the degree of the
polarisation determines the orbit of the polarisation vector and here
in general it does not. This fact is significant for future studies
dealing with classes of polarisations not covered in this paper.

\begin{proposition}\label{hdorbits}
  Let $h_d\in L_{2t}$ be primitive of length $2d>0$ and $\div(h_d)=f$. 
  We put
\begin{align*}
g&=(\sfrac{2t}f, \sfrac{2d}f),&  w&=(g,f),&  g&=wg_1,& f&=wf_1. 
\end{align*}
Then
\[
2t=fgt_1=w^2f_1g_1t_1\text{\qquad and \qquad}
2d=fgd_1=w^2f_1g_1d_1
\]
where $(t_1,d_1)=(f_1,g_1)=1$. 
\begin{itemize}
\item[{(i)}] If $g_1$ is {\bf even}, then such an $h_d$ exists if and
  only if $(d_1,f_1)=(f_1,t_1)=1$ and $-d_1/t_1$ is a quadratic
  residue modulo $f_1$.  Moreover the number of
  $\Tilde\Orth(L_{2t})$-orbits of $h_d$ with fixed $f$ (if at least one
  $h_d$ exists) is equal to
\[
w_+(f_1)\phi(w_-(f_1))\cdot 2^{\rho(f_1)},
\]
where $w=w_+(f_1)w_-(f_1)$ and $w_+(f_1)$ is the product of all powers
of primes dividing $(w,f_1)$, $\rho(n)$ is the number of prime factors
of $n$ and $\phi(n)$ is the Euler function. 
\item[{(ii)}] If $g_1$ is {\bf odd}, and $f_1$ is {\bf even} or $f_1$
  and $d_1$ are {\bf both odd}, then such an $h_d$ exists if and only
  if $(d_1,f_1)=(t_1, 2f_1)=1$ and $-d_1/t_1$ is a quadratic residue
  modulo $2f_1$.  The number of $\Tilde\Orth(L_{2t})$-orbits of such
  $h_d$ is equal to
\begin{alignat*}{2}
w_+(f_1)\phi(w_-(f_1))\cdot
2^{\rho(f_1/2)}&\qquad &\text{if }f_1&\equiv 0 \mod 2\\
\intertext{and to}
w_+(f_1)\phi(w_-(f_1))\cdot
2^{\rho(f_1)}&\qquad &\text{if }f_1&\equiv d_1\equiv 1\mod 2. 
\end{alignat*}
\item[{(iii)}] If $g_1$ and $f_1$ are {\bf both odd} and $d_1$ is
  {\bf even}, then such an $h_d$ exists if and only if $(d_1,f_1)=(t_1,
  2f_1)=1$, $-d_1/(4t_1)$ is a quadratic residue modulo $f_1$ and $w$
  is odd.  The number of $\Tilde\Orth(L_{2t})$-orbits of such $h_d$ is
  equal to
\[
w_+(f_1)\phi(w_-(f_1))\cdot 2^{\rho(f_1)}. 
\]
\item[{(iv)}] For $c$ a suitable integer, determined mod~$f$ and
  satisfying $(c,f)=1$, and $b=(d+c^2t)/f^2$, we have
\[
(h_d)^\perp_{L_{2t}}\cong 2U\oplus 2E_8(-1)\oplus B
\quad\text{with}\quad
B=
\begin{pmatrix}
-2b&c\frac{2t}f\\
c\frac{2t}f&-2t
\end{pmatrix}. 
\]
The form $B$ is a negative definite binary quadratic form of
determinant ${4dt}/{f^2}$.  The greatest common divisor of the
elements of $B$ is equal to $g_1(\frac{2b}{g_1},w)$. 
\end{itemize}
\end{proposition}
\begin{proof}
  A primitive vector $h_d$ with $(h_d, L_{2t})=f\ZZ$ acn be written
\[
h_d=fv+cl_t
\]
 where $v\in 3U\oplus 2E_8(-1)$. The coefficient $c$ is coprime
  to $f$ because $h_d$ is primitive.  According to Eichler's
  criterion (Lemma~\ref{Eichler}) the $\Tilde\Orth(L_{2t})$-orbit of
  $h_d$ is uniquely determined by $h_d^*\equiv \frac {c}{f}l_t\mod
  L_{2t}$.  Therefore it is determined by $c$ mod $f$ because the
  discriminant group of $L_{2t}$ is cyclic. 

We put  $v^2=2b$.  Then $2d=2bf^2-2c^2t$, or
\begin{equation}\label{c-eq}
2f_1b=g_1(d_1+c^2t_1)
\end{equation}
with $(f_1,g_1)=(t_1,d_1)=(c,f_1)=1$. If a $c$ coprime to $f$ and
satisfying Equation~\eqref{c-eq} exists, then because
$(f_1,t_1)|g_1d_1$ we must have $(f_1, t_1)=1$, and similarly we must
have $(f_1,d_1)=1$. 

\smallskip
\noindent(i) First we consider the case when $g_1$ is {\bf
even}. Equation~(\ref{c-eq}) is equivalent to the congruence
\begin{alignat}{2}\label{c1-cong0}
t_1c_1^2 &\equiv -d_1\mod f_1,\qquad & c_1&\equiv c \mod f_1. 
\end{alignat}
If $g_1$ is even, then $f_1$ is odd. Since $(t_1,f_1)=(d_1,f_1)=1$
then $-d_1/t_1\mod f_1$ is invertible modulo $f_1$.  If $-d_1/t_1$ is
not a quadratic residue modulo $f_1$ then the
congruence~(\ref{c1-cong0}) has no solutions, so we assume that
$-d_1/t_1$ is a quadratic residue modulo $f_1$.  Because $f_1$ is odd,
the number of solutions $c_1$ of (\ref{c1-cong0}) taken modulo $f_1$
is equal to
\[
\#\{x \mod f_1\mid x^2\equiv 1 \mod f_1\}=2^{\rho(f_1)}. 
\]
Let us calculate the number of solutions $c$ modulo $f$ where
$f=wf_1$.  Any solution $c_1$ is coprime to $f_1$.  Let us put
\[
c=c_1+(x+yw_-(f_1))f_1
\]
where $x$ is taken mod~$w_-(f_1)$ and $y$ is taken mod~$w_+(f_1)$. 

We note that $(f_1, w_{-}(f_1))=1$.  We have $(c,wf_1)=1$ if and only
if $(c_1+xf_1, w_-(f_1))=1$.  For any fixed $c_1$ the numbers
$c_1+xf_1$ form the full residue system modulo $w_-(f_1)$ if $x$ runs
modulo $w_-(f_1)$.  Therefore for any fixed $c_1$ there are exactly
$w_+(f_1)\phi(w_-(f_1))$ solutions $c$ modulo $f=f_1w$ such that
$c\equiv c_1\mod f_1$ and $(c,f)=1$.  This finishes the proof of~(i). 

\smallskip
\noindent (ii) Consider the case when $g_1$ is {\bf odd}. In this case $t_1$
is also odd, since $(g_1,2f_1)=1$ so $d_1+c^2t_1\equiv
  0\mod 2f_1$: so if $t_1$ is even then $d_1$ is also even. 

If $f_1$ is {\bf even}, then because $f_1^2$ is divisible by $2f_1$
Equation~(\ref{c-eq}) is equivalent to the congruence
\begin{alignat}{2}\label{c1-cong1-1}
t_1c_1^2  &\equiv -d_1\mod 2f_1,\qquad&  c_1&\equiv c \mod f_1. 
\end{alignat}
Since $(f_1,d_1)=1$ we know that $d_1$ is odd. Therefore 
\[
\#\{c_1\mod f_1\ |\ c_1^2\equiv -d_1/t_1\mod 2f_1\}=2^{\rho(f_1/2)}
\]
if $-d_1/t_1$ is a quadratic residue mod~$2f_1$. 
The rest is similar to the case~(i). 

If $f_1$ and $d_1$ are {\bf both odd}, then Equation~(\ref{c-eq}) is
equivalent to the congruence
\begin{alignat}{2}\label{c1-cong1-2}
t_1c_1^2  &\equiv -d_1\mod 2f_1,\qquad&  c_1&\equiv c \mod 2f_1. 
\end{alignat}
The residue $-d_1/t_1\mod 2f_1$ is invertible and it is
always $1$ modulo $2$.  Therefore the number of solutions modulo
$2f_1$ is equal to $2^{\rho(f_1)}$ and they are all different modulo
$f_1$.  We put
\[
c=c_1+(x+yw_-(2f_1)2^{\delta_2(w)-1})2f_1
\] 
taking $x$ mod $w_-(2f_1)$ and $y$ mod $w_-(2f_1)2^{\delta_2(w)-1}$,
where $2^{\delta_2(w)}$ is the $2$-factor of $w$ for $w$ even and is
$2$ if $w$ is odd (so the factor $2^{\delta_2(w)-1}$ only appears if
$w$ is even). For even or odd $w$ we obtain the same formula for the
number of solutions $c$. This proves~(ii). 
\smallskip

\noindent (iii) Consider the case when $g_1$ and $f_1$ are {\bf both odd}, and
$d_1$ is {\bf even}. Equation~(\ref{c-eq}) is equivalent to the
congruence~(\ref{c1-cong1-2}) and $c_1$ is always even, i.e.,
$c_1=2c_2$ and
\[
c_2^2\equiv -(d_1/2)/(2t_1)\mod f_1
\]
and $c_2$ is considered modulo $f_1$.  In particular $w$ must be
odd since otherwise $2|(c,w)$.  For odd $w$ we choose
\[
c=2c_2+(x+w_-(f_1)y)2f_1
\]
with $x$ taken mod $w_-(f_1)$ and $y$ taken $w_-(f_1)$, and this
completes the proof of~(iii). 
\smallskip

\noindent (iv) Fix a representative of the $\Tilde\Orth(L_{2t})$-orbit of $h_d$
of the form
\begin{equation}\label{hd}
h_d=fe_1+fbe_2+cl_t\in U\oplus \latt{-2t}
\end{equation}
where $e_1$ and $e_2$ form a usual basis of the hyperbolic plane $U$
($e_1^2=e_2^2=0$ and $e_1\cdot e_2=1$).  The orthogonal complement of
$h_d$ in $U\oplus \latt{-2t}$ is a lattice $L_B$ of rank $2$
\begin{equation}\label{LB}
L_B=(h_d)_{U+\latt{-2t}}^\perp
=\latt{e_1-be_2,\  c\sfrac{2t}{f} e_2+l_t}. 
\end{equation}
with the quadratic form $B$ as in~(iv). Both vectors are orthogonal to
$h_d$: they form a basis because using them one can reduce to zero
the coordinates at $e_1$ and at $l_t$.  In the notation above we
obtain
\[
B=
g_1\begin{pmatrix} -2b/g_1& cwt_1\\
cwt_1& -w^2f_1t_1
\end{pmatrix}. 
\]
We have $(cwt_1, w^2f_1t_1)=wt_1(c,wf_1)=wt_1$.  The greatest common
divisor of the elements of $B$ is equal to $g_1(\frac{2b}{g_1}, w)$
because $2b/g_1$ and $t_1$ are coprime. 
\end{proof}
\begin{corollary}
Let us assume that $w=1$. If there exists a primitive
vector $h_d\in L_{2t}$ such that $h_d^2=2d$ and $\div(h_d)=f$,
then all such vectors belong to the same $\Orth(L_{2t})$-orbit. 
\end{corollary}
\begin{proof}
  The natural projection $\Orth(L_{2t})\to\Orth(D(L_{2t}))$ is surjective
  (see~\cite{Nik}). Furthermore (see~\cite{GH})
\[
\Orth(D(L_{2t}))\imic \{\,x\mod 2t\mid x^2\equiv 1 \mod 4t\,\}. 
\]
Therefore for $w=1$ all solutions $c\mod f$ of the
congruences (\ref{c1-cong0}) and (\ref{c1-cong1-1}) are equivalent
modulo the action of this abelian $2$-group. 
\end{proof}

\begin{example}\label{splitpolarisation}
  Let $f=1$.  From the Proposition~\ref{hdorbits} it follows that for
  any $t$ and $d$ there is only one $\Tilde\Orth(L_{2t})$-orbit of
  primitive vectors $h_d$ with $\div(h_d)=1$.  Moreover $c=0$ and so
\begin{equation}\label{defineLtd}
(h_d)^\perp_{L_{2t}}\cong L_{2t,2d}
=2U\oplus 2E_8(-1)\oplus\latt{-2t}\oplus\latt{-2d}. 
\end{equation}
\end{example}

\begin{definition}\label{def:split}
We call a polarisation determined by a primitive vector $h_d$ {\it split\/}
if $\div(h_d)=1$ and {\it non-split\/} otherwise.
\end{definition}

The name comes from the fact that $\div(h_d)=1$ 
is the only case when the matrix $B$ is
diagonal because $c$ can be zero and coprime to $f$ only if $f=1$. 
The split case will be the main subject of much of the rest of this paper.

\begin{example}\label{f=2}
  Let $f=2$.  In this case $c$ is odd, so we may take $c=1$.  A
  constant $b$ and a vector $h_d$ exist if and only if $d+t\equiv
  0\mod 4$.  Moreover the $\Tilde\Orth(L_{2t})$-orbit of $h_d$ is unique
  because $D(L_{2t})$ is cyclic and thus contains only one element of
  order $2$. 
\end{example}
If $f>2$, then Proposition~\ref{hdorbits} shows that the number of
orbits is zero or strictly greater than one. Thus the cases $f=1$ and
$f=2$ are special in that they are the only cases where the degree determines the 
polarisation uniquely.

\begin{example}\label{(d,t)=1} 
  Let $d$ and $t$ be coprime. Examples~\ref{splitpolarisation} and
  \ref{f=2} give us the full classification of possible $h_d\in
  L_{2t}$ (in particular if $t=1$ or
  $d=1$), since if $(t,d)=1$ then $f=\div(h_d)=1$ or $2$. 
\end{example}

In the next proposition we show that if $w=1$, then the groups
$\Tilde\Orth(L_{2t}, h_d)$ and $\Orth(L_{2t}, h_d)$ have rather clear
structure. 

\begin{proposition}\label{ogroup} Let $h_d\in L_{2t}$ be a primitive vector
  such that $h_d^2=2d$ and $\div(h_d)=f$. Assume that $w=1$,
  {i.e.} $f$ and $(\frac{2t}f, \frac{2d}f)$ are coprime. Then
\begin{itemize}
\item[{(i)}]
$\Tilde\Orth(L_{2t}, h_d)\cong \Tilde\Orth((h_d)^\perp_{L_{2t}})$. 
\item[{(ii)}] The factor group $\Orth(L_{2t}, h_d)/\Tilde\Orth(L_{2t},
  h_d)$ is an abelian $2$-group, which is of order $ 2^{\rho({t}/f)}$ if $f$ is
  odd. If $f$ is even the order is equal to $2^{\rho({2t}/f)+\delta}$,
  where
\[
\delta=
\begin{cases}
\hphantom{-}0&\text{\quad if }(2t/f ) \equiv 1\mod 2 \text{ or }(2t/f)\equiv 4\mod 8,\\
-1&\text{\quad if } (2t/f) \equiv 2\mod 4,\\
\hphantom{-}1&\text{\quad if } (2t/f) \equiv 0\mod 8. 
\end{cases}
\]
\end{itemize}
\end{proposition}
\begin{proof}
  We may take $h_d$ in the form (\ref{hd}).  We can fix a basis
  ${k_1,k_2}$ of $L_B^\vee$ given by
\begin{alignat}{2}
k_1&=\frac{f}{2d}h_d-e_2&&=\frac{f}{2d}\bigl(fe_1+bfe_2+cl_t\bigl)-e_2,\notag\\
k_2&=\frac{c}{2d}h_d+\frac{1}{2t}l_t&&=
\frac{f}{2d}\bigr(ce_1+cbe_2+\frac{bf}tl_t).\notag 
\end{alignat}
Up to sign, this is dual to the basis fixed in~\eqref{LB}.  We put
\[
k_3=fk_2-ck_1=ce_2+\frac{f}{2t}l_t. 
\]
If $v\in L^\vee$ we shall denote by $\bar v$ the corresponding element
in the discriminant group $D(L)=L^\vee/L$.  We note that the orders of
$\bar{k}_1$ and of $\bar{k}_3$ in $D(L_B)=D((h_d)^\perp_{L_{2t}})$ (see
\eqref{LB}) are equal to $\frac{2d}f$ and $\frac{2t}f$ respectively. 
Moreover $\bar{k}_1\cdot \bar{k}_3=0$.  Let us calculate the order of
the intersection of the subgroups generated by $\bar{k}_1$ and
$\bar{k}_3$ in $D(L_B)$. If $n\bar{k}_1\in \latt{\bar{k}_3}$ then
$n=g_1d_1n_1$ because $m\bar{k}_3$ does not contain the
$e_1$-component.  Therefore $n\bar{k}_1\equiv
(\frac{n_1c}{w}l_t+xe_2)\mod L_B$ where $x\in \ZZ$ (see \eqref{LB}),
and $|\latt{\bar{k}_1}\cap \latt{\bar{k}_3}|=w$.  It follows that
$\bar{k}_1$ and $\bar{k}_3$ form a basis of $D(L_B)$ if $w=1$. 

As in the beginning of the section we consider the following series of
lattices
\[
\latt{ h_d} \oplus h_d^\perp < L_{2t} < L_{2t}^\vee
< \latt{ h_d^\vee}\oplus (h_d^\perp)^\vee,
\]
where $h_d^\vee=\frac{1}{2d}h_d$ and $h_d^\perp=(h_d)^\perp_{L_{2t}}\cong
2U\oplus 2E_8(-1)\oplus L_B$.  The subgroup $H=L_{2t}/(\latt{ h_d}
\oplus h_d^\perp))<D(h_d)\oplus D(L_B)$ has order $\frac{2d}f$. It is
generated by the element $\bar{k}_1-f\bar h_d/2d$. Therefore the
projection
\[
p(H)=p_{h_d^\perp}(H)=\latt{ {\bar k}_1}
\]
is the subgroup generated by ${\bar k}_1$. It follows if $w=1$ the
discriminant group is
\[
D(h_d^\perp)=\latt{ \bar{k}_1} \oplus \latt{ \bar{k}_3}=
p(H)\oplus \latt{ \bar{k}_3}. 
\]
According to Lemma~\ref{ols}
\[
\Orth(L_{2t}, h_d)\cong
\{\gamma\in \Orth(h_d^\perp)\mid \bar \gamma|_{p(H)}=\id\}. 
\]
Let us consider an element $\gamma\in \Orth(h_d^\perp)$ satisfying
$\bar \gamma|_{p(H)}=\id$ as an element of $\Orth(L_{2t}, h_d)$ (i.e., we
put $\gamma(h_d)=h_d$).  According to the decomposition above
$\gamma\in \Tilde\Orth(L_{2t}, h_d)$ if and only if $\bar\gamma(\bar
{k}_3)=\bar {k}_3$.  Therefore $\Tilde \Orth(L_{2t}, h_d)\cong \Tilde
\Orth(h_d^\perp)$. 

We note that the natural projection $\Orth(h_d^\perp)\to
\Orth(D(L_B))$ is surjective (see~\cite{Nik}).  Therefore according to
Lemma~\ref{ols}
\[
\Orth(L_{2t}, h_d)/\Tilde\Orth(L_{2t}, h_d) \cong\{\gamma\in
\Orth(h_d^\perp)\mid\bar \gamma|_{p(H)}=\id\}/
\Tilde\Orth(h_d^\perp)\cong \Orth(\latt{ \bar k_3}),
\]
where $\latt{ \bar k_3}=\{n\bar k_3\mid n\mod \frac{2t}f\}$ and
$\bar k_3^2\equiv -\frac{f^2}{2t}\mod 2$.  Therefore
\begin{align*}
\Orth(\latt{ \bar k_3})&\cong
\{x\mod \sfrac{2t}f\mid x^2\bar k_3^2\equiv  \bar k_3^2\mod 2\}\\ 
&=\{x \mod \sfrac{2t}f\mid
x^2f\equiv f \mod 2\sfrac{2t}f\} 
\end{align*}
We supposed that $w=1$. Therefore
$f=f_1$ and $g=g_1$ are coprime.  We have $\frac{2t}f=g_1t_1$ with
$(f_1,t_1)=1$ (see Proposition~\ref{hdorbits}).  It follows that the
group $\Orth(\latt{ \bar k_3})$ is isomorphic to the group
\[
\{x \mod \sfrac{2t}f\mid x^2\equiv 1 \mod 2^{\varepsilon(f)}\sfrac{2t}f\},
\]
where $\varepsilon(f)=1$ if $f$ is odd (in this case $\frac{2t}f$ is
even) and $\varepsilon(f)=0$ is $f$ is even. The last group is well-known (compare
with~\cite{GH}). 
\end{proof}
\begin{corollary}\label{o=tilde}
  We have that $\Orth(L_{2t}, h_d)\cong \Tilde\Orth(L_{2t},h_d)$ in the
  following three cases: $f$ is odd and $f=t$; or $f=2t$; or $f=t$ and
  $2d/f$ is odd. 
\end{corollary}
\begin{proof} If $f=t$ or $f=2t$, then $g=(\frac{2t}f, \frac{2d}f)=1$
  or $2$.  If $g=1$, then $w=1$. If $g=2$ then $w=(f,g)=1$ for odd
  $f$ and for even $f$ such that $(2d)/f$ is odd. In all these case
  the index $[\Orth(L_{2t}, h_d): \Tilde\Orth(L_{2t},h_d)]=1$ according
  Proposition~\ref{ogroup}. 
\end{proof}

\begin{remark}
The condition $w=1$, i.e., that $f$ and $(\frac{2t}f,
\frac{2d}f)$ are coprime, is valid for any $f$ if $(2t,2d)$ is square
free. In particular this condition is true for any vector $h_d$ if
$2t$ is square free.  
\end{remark}

\begin{remark}
  The finite group $\Orth(D((h_d)^{\perp}_{L_{2t}}))$ is cyclic for
  any $h_d$ with $\div (h_d)=f$ if $g_1=(\frac{2t}f, \frac{2d}f)=1$. 
  If $g_1>1$ the discriminant group is not cyclic, but it is the
  orthogonal sum of two cyclic groups if $w=1$. 
  Proposition~\ref{ogroup} shows that we can consider this rather
  general case as a regular one. 
\end{remark}

Since the classification of polarisation types in this section depends
only on the discriminant group it immediately gives an identical
classification for polarisations of deformations of generalised Kummer
varieties. 

\section{Modular forms and root systems}\label{sec:modular}

For the rest of the paper we restrict to a special class of symplectic
$4$-folds. We consider the case of deformation $\Kthree^{[2]}$
manifolds with polarisation of degree $2d$ of split type, as in
Example~\ref{splitpolarisation} above. We denote an irreducible
component of the corresponding moduli space by
$\cM^{[2],\text{split}}_{2d}$. (We do not know whether there is only
one irreducible component: see Question~\ref{qu:irreducible}.) 
According to Theorem~\ref{theo:cover} we have a dominant map 
\[
\cM^{[2],\text{split}}_{2d} \To \Tilde\Orth(L_2,h_d) \backslash
\Omega_{h_d}. 
\]
In this case, where $t=1$, we have 
\[
\Tilde\Orth(L_{2,2d})=\Tilde\Orth(L_2,h_d)=\Orth(L_2,h_d)
\]
by Proposition~\ref{ogroup}(i) and Corollary~\ref{o=tilde}, where
$L_{2,2d}$ is as defined in Equation~\eqref{defineLtd}. In
particular the vertical map in Theorem~\ref{theo:cover} is of
degree~$1$, and an affirmative answer to Question~\ref{qu:irreducible}
would imply that global Torelli holds for deformation $\Kthree^{[2]}$
manifolds with polarisation of split type. 

It is more convenient to express this quotient in terms of the
symmetric domain $\cD(L_{2,2d})$ defined in Equation~\eqref{cD} above. 
Recall that $\Tilde\Orth^+(L_{2,2d})$ is the index~$2$ subgroup of
$\Tilde\Orth(L_{2,2d})$ that preserves $\cD(L_{2,2d})$. Then
\[
\Tilde\Orth(L_2,h_d) \backslash
\Omega_{h_d}=\Tilde\Orth^+(L_{2,2d})\backslash \cD(L_{2,2d}). 
\]
In the rest of this paper we study the Kodaira dimension of the
locally symmetric variety 
$\Tilde\Orth^+(L_{2,2d})\backslash \cD(L_{2,2d})$ and of the moduli
space $\cM^{[2],\text{split}}_{2d}$. 
 
\begin{theorem}\label{gentype}
  The variety $\cM^{[2],\text{split}}_{2d}$ is of general type if $d\ge
  12$.  Moreover its Kodaira dimension is non-negative if $d=9$ and
  $d=11$. 
\end{theorem}

In the case $d=11$ the split and non-split cases behave very
differently. Debarre and Voisin~\cite{DV} give a geometric
construction of a $20$-dimensional family of irreducible symplectic
manifolds with non-split polarisation with $d=11$, and from their
construction it follows immediately that the moduli space is unirational.

Let $L$ be an even integral lattice of signature $(2,n)$ with $n\ge
3$.  A modular form of weight $k$ and character $\chi\colon
\Gamma\to\CC^*$ for a subgroup $\Gamma<\Orth^+(L)$ of finite index is
a holomorphic function $F\colon\cD_L^\bullet\to \CC$ on the affine
cone $\cD_L^\bullet$ over $\cD_L$ such that
\begin{equation*}\label{mod-form}
  F(tZ)=t^{-k}F(Z)\quad \forall\,t\in \CC^*\text{\qquad and\qquad}
  F(gZ)=\chi(g)F(Z)\quad \forall\,g\in \Gamma. 
\end{equation*}
A modular form is a cusp form if it vanishes at
every cusp. For applications, we require the order of vanishing to be
at least~$1$ (both here and in~\cite{GHS:K3}, although it is not
stated explicitly there). In general this is a slightly stronger
requirement because the order of vanishing might be a rational number
less than~$1$. However, it is easy to check that for trivial character and
character~$\det$, which are the cases used here and in~\cite{GHS:K3},
the vanishing order at any cusp is an integer. 

We denote the linear spaces of modular and cusp forms of weight $k$
and character $\chi$ for $\Gamma$ by $M_k(\Gamma,\chi)$ and
$S_k(\Gamma,\chi)$ respectively. 

The next theorem follows from  the results obtained in \cite{GHS:K3}. 
\begin{theorem}\label{modvar}
  Suppose there exists a non-zero cusp form $F_a$ of some weight
  $a<20$ and character $\det$ with respect to the modular group
  $\Tilde\Orth^+(L_{2,2d})$.  Then the modular variety
  $\cM^{[2],\text{split}}_{2d}$ is of general type. 

  If there exists a non-zero cusp form $F_{20}$ of weight $20$ and
  character $\det$ then $\cM^{[2],\text{split}}_{2d}$ has non-negative
  Kodaira dimension. 
\end{theorem}
\begin{proof} 
  $\cM^{[2],\text{split}}_{2d}$ is a quasi-projective variety of
  dimension $20$. It has a toroidal compactification having only
  canonical singularities, by \cite[Theorem 2]{GHS:K3}.  By
  \cite[Theorem 1.1]{GHS:K3}, the variety
  $\cM^{[2],\text{split}}_{2d}$ is of general type if there exists a
  non-zero cusp form $F_a\in S_a(\Tilde\Orth^+(L_{2,2d}))$ of weight
  $a<20$ that vanishes along the ramification divisor of the
  projection
\[
\pi\colon \cD(L_{2,2d}) \To \Tilde\Orth^+(L_{2,2d})\backslash \cD(L_{2,2d}). 
\]
We note that according to \cite[Corollary 2.13]{GHS:K3} the
ramification divisor is determined by the elements $\sigma\in
\Tilde\Orth(L_{2,2d})$ such that $\sigma$ or $-\sigma$ is a reflection
with respect to a vector $r\in L_{2,2d}$.  We classified those
reflections using the results of \cite[\S 3]{GHS:K3}. 

Let $F_a\in S_a(\Tilde\Orth^+(L_{2,2d}),\det)$ be of weight $a<20$ and
suppose that $\sigma\in\Tilde\Orth(L_{2,2d})$ defines a component of the
ramification divisor.  Then
\[
F_a(\pm\sigma(Z))=\det(\pm \sigma)\cdot F_a(Z)=-F_a(Z)
\]
because $\det(-\sigma)=(-1)^{20}\det(\sigma)=-1$.  Therefore the cusp
form $F_a$ with character $\det$ automatically vanishes on the
ramification divisor. 

If $a=20$ and $F_{20}$ vanishes along the ramification divisor
of~$\pi$ then $F_{20}$ determines a section of the canonical bundle by
a well-known result of Frei\-tag~\cite[Hilfssatz 2.1, Kap. III]{Fr}. 
\end{proof}

One can estimate the obstructions to continuing of the pluricanonical
forms across the ramification divisor using the exact formula for
Mumford--Hirzebruch volume of the corresponding orthogonal groups (see
\cite{GHS:Vol}). But this approach only gives good results for locally
symmetric varieties of orthogonal type if the dimension is quite
large: at least $33$ in the cases considered in \cite{GHS:Orth}. 

If the dimension of the modular variety is smaller than $26$ we can
use the quasi pull-back (see Equation~\eqref{pb} below) of the
Borcherds modular form $\Phi_{12}$ to construct cusp forms of small
weight. The Borcherds form is a modular form of weight~$12$ for
$\Orth^+(II_{2,26})$, where $II_{2,26}$ is the unimodular lattice
$2U\oplus 3E_8(-1)$. 

$\Phi_{12}(Z)=0$ if and only if there exists $r\in
II_{2,26}$ with $r^2=-2$ such that $(r,Z)=0$. Moreover, the
multiplicity of the rational quadratic divisor in the divisor of zeros
of $\Phi_{12}$ is~$1$ (see \cite{Bo}).  This form generates very
important functions on the moduli spaces of polarised $\Kthree$
surfaces (see \cite{BKPS}, \cite{Ko} and \cite{GHS:K3}).  In the context
of the moduli space of symplectic manifolds we can use the following
specialisation of the quasi pull-back. 

The Weyl group of $E_8$ acts transitively on the roots of $E_8$.  If
$v$ is a root of $E_8(-1)$ then $v^\perp_{E_8(-1)}\cong E_7(-1)$.  Let
$l\in E_7(-1)$ satisfy $l^2=-2d$.  The choice of $v$ and $l$
determines an embedding of $L_{2,2d}$ into $II_{2,26}$.  The embedding
of the lattice also gives us an embedding of the domain
$\cD(L_{2,2d})\subset \PP (L_{2,2d}\otimes \CC)$ into $\cD(II_{2,26})
\subset \PP(II_{2,26} \otimes \CC)$. 

We put $R_l=\{r\in E_7(-1)\mid r^2=-2,\ (r, l)=0\}$, and $N_l=\# R_l$. 
(It is clear that $N_l$ is even.)  We note that $R_l$ is the set of
roots orthogonal to the sublattice $\latt{v}\oplus\latt{l}$ in
$E_8(-1)$.  Then the quasi pull-back of $\Phi_{12}$ is given by the
following formula (see \cite{BKPS}):
\begin{equation}\label{pb}
  \left. F_l=
  \frac{\Phi_{12}(Z)}{
  \prod_{\{ \pm r\}\in R_l} (Z, r)}
  \ \right\vert_{\cD_{L_{2,2d}}}
  \in M_{12+\frac{N_l}2}(\Tilde\Orth^+(L_{2,2d}),\, \det). 
\end{equation}
It is a non-trivial modular form of weight $12+\frac{N_l}2$. By
\cite[Theorem 6.2]{GHS:K3} it is a cusp form if $N_l$ is non empty. In
\cite{GHS:K3} we proved this for $l\in E_8(-1)$.  But in the proof we
used only the fact that any isotropic subgroup of the discriminant
form of the lattice $2U\oplus 2E_8(-1)\oplus\latt{-2d}$ is
cyclic (see \cite[Theorem 4.2]{GHS:K3}).  The same is true for
$L_{2,2d}$ because its discriminant group is cyclic (see \S 2).  The
weight of $F_l$ is smaller than $20$ if $N_l<16$. 

The problem therefore is to determine the $d$ for which such a vector
exists. Sufficient conditions are given in Theorem~\ref{mainineq}
below. We apply the method used in the proof of 
\cite[Theorem 7.1]{GHS:K3}. We first need some properties of the lattice $E_7$. 

\begin{lemma}\label{E7} 
  The Weyl group $W(E_7)$ acts transitively on the sets of sublattices
  of $E_7$ of types $A_1\oplus A_1$ or $A_2$. 
\end{lemma}
\begin{proof}
  $W(E_7)$ acts transitively on the roots. Moreover
  $(A_1)^\perp_{E_7}\cong D_6$ and $W(D_6)$ acts transitively on its
  roots. This proves the $A_1\oplus A_1$ case. 

  Let $A_2^{(1)}$, $A_2^{(2)}$ be two different copies of $A_2$ in
  $E_7$.  Without loss of generality we can assume that they have a
  common root $a$, i.e. $A_2^{(1)}=A_2(a,c)=\ZZ a+\ZZ c$ and
  $A_2^{(2)}=A_2(a,d)$, where $a\cdot c=a\cdot d=-1$.  Any
  $A_2$-lattice contains six roots
\[
R(A_2(a,c))=\{\,\pm a,\ \pm c,\ \pm(a+c)\,\}
\]
and it is generated by any pair of linearly independent roots.  If
$c\cdot d=-1$ then $(a+d)\cdot c=-2$, $c=-(a+d)$ and $A_2(a,c)=
A_2(a,d)$. Therefore $c\cdot d=0$ or $1$.  (We recall that for any two
non-collinear roots $u$ and $v$ one has $|u\cdot v|\le 1$.) 

If $c\cdot d=1$ then $(c-d)^2=2$ and $(c-d)\cdot a=0$.  The reflection
$\sigma_{c-d}$with respect to the root $(c-d)$ transforms $A_2^{(1)}$ into
$A_2^{(2)}$:
\[
\sigma_{(c-d)}(c)=c-(c\cdot (c-d))(c-d)=d,\qquad
\sigma_{(c-d)}(a)=a. 
\]

If $c\cdot d=0$ then $A_2(a,c)+A_2(a,d)$ is a root lattice of type
$A_3$ with $12$ roots:
\[
R(A_3(a,c,d))
=\pm \{a,\  c,\  d,\  a+c,\  a+d,\  a+c+d\}. 
\] 
The roots $\pm(a+c+d)$ are the only roots in $A_3(a,c,d)$
orthogonal to $a$.  We have $\sigma_{a+c+d}(c)=-(a+d)$.  To find a
reflection $\sigma$ such that $\sigma(a)=a$ and $\sigma(-(a+d))=d$ we
have to go outside of $A_3=A_2(a,c)+A_2(a,d)$.  We recall that $E_7$
contains $126$ roots (see \cite{Bou}):
\begin{align*}
R(E_7)=&\{\,\pm(e_7-e_8)\,\}\cup\{\,\pm e_i\pm e_j\mid 1\le i< j \le
6\,\}\\ 
&\cup\{\,\pm \sfrac{1}{2}(e_7-e_8)+\sfrac{1}{2} \sum _{i=1}^6
(-1)^{\nu(i)} e_i\mid\sum_{i=1}^6 \nu(i) \text{ is even}\,\}
\end{align*}
where $e_i$ form the usual Euclidian basis in $\ZZ^8$.  Without loss
of generality we can assume that $a=e_7-e_8$.  Since $a\cdot d=-1$, we
see that
\[
d=-\sfrac1{2}(e_7-e_8)+\sfrac{1}2 \sum _{i=1}^6 (-1)^{\nu(i)} e_i. 
\]
We put
$r_j=(-1)^{\nu(j)}e_j+(-1)^{\nu(j+1)}e_{j+1}$ for $j=1$, $3$ and $5$. 
We obtain
\[
\sigma_{r_5}\circ \sigma_{r_3}\circ \sigma_{r_1}(d)=d-r_1-r_3-r_5=-(a+d). 
\]
\end{proof}
\begin{corollary}
We have
\begin{itemize}
\item[(i)] The Weyl group $W(E_8)$ acts transitively on all its
  sublattices of types $3A_1$ and $A_1\oplus A_2$ in $E_8$. 
\item[(ii)] The class number of the lattices $A_5$ and $A_1\oplus D_4$
  is equal to one. 
\item[(iii)] The sublattices $4A_1$ in $E_8$ form two orbits with respect
  to $W(E_8)$. 
\end{itemize}
\end{corollary}
\begin{proof}
  (i) follows from the fact that $(A_1)^\perp_{E_8}\cong
  E_7$.  To prove (ii) we note that $A_1\oplus D_4$ is a
  maximal lattice because its discriminant group does not contain any
  isotropic vectors.  (The square of any element in $D_4^\vee/D_4$ is
  equal to $1/2$ modulo $2$.) Furthermore
\begin{align}\label{orthcompl}
(A_1)^\perp_{E_7}&\cong D_6,& (A_1)^\perp_{D_6}&\cong A_1\oplus D_4,
&(A_2)^\perp_{E_7}&\cong A_5. 
\end{align}
To see these one has to use the extended Dynkin diagram of the
corresponding root lattice and to take into account the maximality of
$D_6$, $A_5$ and $A_1\oplus D_4$.  The discriminant quadratic form is
the invariant of the genus of an even quadratic lattice. Therefore if
$M$ is a lattice in the genus of $A_1\oplus D_4$ then we can
consider $M\oplus 2A_1$ as a sublattice of $E_7$.  All such $M$ are
isomorphic according to the lemma. The same argument works for $A_5$. 

To prove (iii) we remark that $(3A_1)^{\perp}_{E_8}\cong A_1\oplus D_4$. 
The last lattice contains two orbits of roots. 
\end{proof}

\begin{theorem}\label{mainineq}
  There exists a vector $l$ in $E_7$ of length $2d$ orthogonal to
  at least $2$ and at most $14$ roots if
\begin{equation}\label{mineq1}
30N_{A_1\oplus D_4}(2d)+16N_{A_5}(2d) < 5N_{D_6}(2d)
\end{equation}
or to at least $2$ and at most $16$ roots if
\begin{equation}\label{mineq2}
30N_{A_1\oplus D_4}(2d)+16N_{A_5}(2d) < 6N_{D_6}(2d),
\end{equation}
where $N_L(2d)$ denotes the number of representations of $2d$ by the
lattice $L$. 
\end{theorem}

\begin{proof} 
  Suppose that any vector $l\in E_7$ of length $2d$ is orthogonal to at
  least $16$ roots if it is orthogonal to any.  Let us fix a root $a$
  in $E_7$ orthogonal to $l$.  Therefore $l\in
  a^\perp_{E_7}=D_6^{(a)}$.  The others roots orthogonal to $l$ are
  some roots in $D_6^{(a)}$ ($60$ roots) or roots in $R(E_7)\setminus
  R(D_6^{(a)})$ ($66$ roots).  The last $66$ roots form a bouquet
  $Q_a(16A_2)$ of $16$ copies $A_2(a)$ of $A_2$ centred in $\pm a$. 
  If $l$ is orthogonal to any root from $A_2(a)$ different from $a$,
  then $l$ is orthogonal to the whole lattice $A_2(a)$ and $l\in
  (A_2)_{E_7}^\perp\cong A_5$.  If $l$ is orthogonal to a root in
  $D_6^{(a)}$ then $l\in (2A_1)_{E_7}^\perp\cong A_1\oplus D_4$. 
  Therefore we have
\begin{equation}\label{union}
l\in \bigcup_{i=1}^{30} (A_1\oplus D_4)^{(i)} \cup
\bigcup_{j=1}^{16} A_5^{(j)}. 
\end{equation}
Denote by $n(l)$ the number of components in \eqref{union} containing
the vector $l$.  We have calculated this vector exactly $n(l)$ times
in the sum
\[
30N_{A_1\oplus D_4}(2d)+16N_{A_5}(2d). 
\]
We need to estimate $n(l)$.  We shall consider several cases.  

\smallskip
1. Let $l\cdot c\ne 0$ for any $c\in Q_a(16A_2)\setminus \{\pm a\}$. 
Then $l$ is orthogonal to at least $7$ copies of $A_1$ in $D_6^{(a)}$
and $n(l)\ge 7$.  

\smallskip

Now we suppose that there exist $c\in Q_a(16A_2)\setminus \{\pm a\}$
such that $l\cdot c=0$. Then $l$ is orthogonal to $A_2(a,c)$ which is
one of the $16$ subsystems of the bouquet $Q_a$.  \smallskip

2. If $l$ is orthogonal to only one copy of $A_2$ in $Q_a$, $A_2^{(i)}$ ($6$ roots) then $l$
is orthogonal to at least $5$ copies of $A_1$ in $D_6^{(a)}$.  Thus
$n(l)\ge 6$.  \smallskip

3. If $l$ is orthogonal to exactly two copies of $A_2$ in $Q_a$,
$A_2^{(i)}$ and $A_2^{(j)}$, then $l$ is orthogonal to
$A_3=A_2^{(i)}+A_2^{(j)}$ having $12$ roots.  Thus $l$ is orthogonal
to another $2A_1$ in $D_6^{(a)}$.  But $A_3$ contains one more copy
$A_1$ from $D_6^{(a)}$ orthogonal to $a$ (see the proof of
Lemma~\ref{E7}). Therefore $n(l)\ge 5$. 
\smallskip

4. If $l$ is orthogonal to three or more $A_2^{(i)}$ then their sum
contains at least three $A_1< D_6^{(a)}$ and $n(l)\ge 6$. 
\smallskip

Therefore we have proved that if any $l\in E_7$ with $l^2=2d$ is
orthogonal to at least $16$ roots then $n(l)\ge 5$ and
\[
30N_{A_1\oplus D_4}(2d)+16N_{A_5}(2d)\ge 5N_{D_6}(2d). 
\]
If we replace $16$ roots by $18$ roots in the last condition then we
obtain the second inequality of Theorem~\ref{mainineq}
\[
30N_{A_1\oplus D_4}(2d)+16N_{A_5}(2d)\ge 6N_{D_6}(2d). 
\]
\end{proof}

\section{Representations by quadratic forms of odd rank}\label{sec:odd}

To estimate the values of $d$ for which the inequality of Theorem~\ref{mainineq}
is true we need exact formulae for the numbers $N_{A_1\oplus D_4}(2d)$
and $N_{A_5}(2d)$. 

Let $A$ be a symmetric even integral positive definite $m\times
m$ matrix of determinant $\det A=|A|$, and
\[
S(X)=\frac{1}2 A[X]= \frac{1}2 {}^tX A X
\] 
be the corresponding quadratic form taking integral values on $\ZZ^m$. 
The genus $\gen S$ of $S$ contains a finite number of classes $S_i$.  The
integral orthogonal group $\Orth(S_i)$ is finite of order
$|\Orth(S_i)|$. One defines the {\it mass\/} of the genus by
\[
\mass (S)=\sum_{S_i\in \gen S}\,|\Orth(S_i)|^{-1}
\]
and the {\it weight\/} of the class $S_i$ in the genus of $S$ by
\[
w(S_i)=|\Orth(S_i)|^{-1}/\mass (S). 
\]
Siegel's main theorem on quadratic forms (see \cite{Si}) tells
us that the number of representations of $t$ by the genus of $S$,
defined by
\[
r(t, \gen S)=\sum_{S_i\in \gen S}\,w(S_i)r(t,S_i)
\]
where $r(t,S_i)$ is the number of the representation of $t$ by the
quadratic form $S_i$, can be written in terms of the local densities
$\alpha_p(t,S)$
\[
r(t, \gen S)=\varepsilon_m\prod_{p\le \infty}\,\alpha_p(t,S),
\]
where $\varepsilon_m=1$ for all $m\ge 2$ except $\varepsilon_2= 1/2$. 
The local densities (or the local measures of the representations of
$t$ by $S$) are defined as follows
\[
\alpha_p(t,S)=\lim_{a\to \infty}
p^{-a(m-1)}\#\{\,X\in (\ZZ/p^a\ZZ)^m\mid S(X)\equiv t \mod p^a\,\}
\]
if $p$ is a finite prime and
\[
\alpha_\infty(t,S)=\lim_{V\to t}\frac{\vol S^{-1}(V)}{\vol V}
=(2\pi)^{\frac{m}2}\Gamma(\sfrac{m}2)^{-1}t^{\frac{m}2-1}|A|^{-\frac{1}2},
\]
where $|A|=\det(A)$, $V$ is a real neighbourhood of $t$ and $\vol$ is
the usual Euclidian volume in $\RR$ or $\RR^m$ (see
\cite[Hilfssatz 26 and (71)]{Si}). 

If the genus of $S$ contains only one class then the Siegel formula
gives us an exact formula for the number $r(t,S)$ of representations
of $t$ by $S$. In his first paper \cite{Si} on the analytic theory of
quadratic forms Siegel found exact formulae for the local densities if
the prime $p$ is not a divisor of the determinant of $A$. If the rank
$m\ge 4$ is even we have the following formula (see
\cite[(11.74)]{Iw})
\begin{equation}\label{evencase}
r(t, \gen S)=a_\infty(t,S)
L(\frac{m}2,\,\chi_{4D})^{-1}\bigl(\sum_{a|t}\chi_{4D}(a)a^{1-\frac{m}2}\bigr)
\cdot \prod_{p|2D}\alpha_p(t,S)
\end{equation}
where $D=(-1)^{m/2}|A|$ is the discriminant of $A$ and
$\chi_{4D}(a)=\left(\frac{4D}{a}\right)$ is the quadratic character. 

Usually the exact computation of the local densities for odd rank $m$
is said to be more complicated: see for example the introduction
to~\cite{Sh}. Here we give a well-organised formula for $r(t, \gen
S)$ for odd rank $m$. For this purpose we use the Zagier $L$-function
$L(s,\Delta)$ and the H.~Cohen numbers $H(n,\Delta)$ (see \cite{C} and
\cite{Za}). In these terms, surprisingly, the exact formula for odd
rank is simpler that the formula~\eqref{evencase} for even rank. 

If $\Delta\equiv 0,1\mod 4$ then $\Delta=Df^2$, where $D$ is the
discriminant of the quadratic field $\QQ(\sqrt{\Delta})$. By
definition (see \cite[(7) and Proposition 3]{Za}) one has
\begin{equation}\label{LsD}
L(s,\Delta)=\frac{\zeta(2s)}{\zeta(s)}
\sum_{n=1}^{\infty}b_n(\Delta)n^{-s},
\end{equation}
where 
$b_n(\Delta)=\#\{x\mod 2n\mid x^2\equiv \Delta \mod 4n\}$,
and 
\[
L(s,\Delta)=
L(s,\chi_D)\sum_{a|f}\mu(a)\left(\frac D{a}\right) a^{-s}
\sigma_{1-2s}\bigl(\sfrac{f}a\bigr),
\]
where $\sigma_s(t)=\sum_{d|t}d^s$ and $\mu(a)$ is the M\"obius
function.  The main advantage of $L(s,\Delta)$ is the fact that it
satisfies a simple functional equation. 

The function
\[
L^*(s,\Delta)=
\begin{cases}
\pi^{-\frac{s}2}\Gamma(\frac{s}2)\Delta^{\frac{s}2}L(s,\Delta)&
\text{if }\Delta>0\\
\pi^{-\frac{s}2}\Gamma(\frac{s+1}2)|\Delta|^{\frac{s}2}L(s,\Delta)&
\text{if }\Delta<0
\end{cases}
\]
has a meromorphic continuation to the whole complex plane
and satisfies the functional equation
\[
L^*(s,\Delta)=L^*(1-s,\Delta)
\]
(see \cite[Proposition 3]{Za}). 
Moreover $L(s,\Delta)$ is entire except for a simple pole
(of residue $\frac{1}2$ if $\Delta=0$ and $1$ otherwise)
if $\Delta$ is a square. 
This function is very useful for calculation of Fourier
coefficients of various Eisenstein series (see \cite{C},
\cite{Za} and \cite{GS-P}). 

To formulate our reorganisation of the Siegel formula for odd rank $m$
we introduce some notation. We write
\[
t=t_At_1t_2^2,
\]
where $t_1$ is square free, $(t_1t_2^2, |A|)=1$ and $t_A$ divides some
power of $|A|$.  We put
\[
D=\disc\QQ\left(\sqrt{(-1)^{\frac{m-1}2}2t|A|}\,\right). 
\] 
We note that $D>0$ if $m\equiv 1\mod 4$ and $D<0$ if $m\equiv 3\mod
4$.  The determinant of $A$ is always even if $m$ is odd and $A$ is
even integral. 
\begin{theorem}\label{rtS}
Let $m=2m_1+1$ and $S(X)=\frac{1}2A[X]$. Then we have 
\begin{multline}\label{L-zeta} 
r(t,\gen S)=\\
(2\pi)^{\frac{m}2}\Gamma\big(\frac{m}2\big)^{-1}t^{\frac{m}2-1}|A|^{-\frac{1}2}
L\big(\frac{m-1}2,\,Dt_2^2\big)\zeta(m-1)^{-1}\\
\cdot \prod_{p|\,|A|}
\frac{1-\chi_D(p)p^{\frac{1-m}2}}{1-p^{1-m}}
\alpha_p(t,S)
\end{multline}
and
\begin{multline}\label{HCohen}
r(t,\gen S)=\\
\left(\frac{t_A}{|D_A|}\right)^{m_1-\frac{1}2}
2^{2m_1-\frac{1}2}
|A|^{-\frac{1}2}
\left|\frac{2m_1}{B_{2m_1}}\right|
(-1)^{[m_1/2]}H(m_1, Dt_2^2)\\
\cdot \prod_{p|\,|A|}
\frac{1-\chi_D(p)p^{\frac{1-m}2}}{1-p^{1-m}}
\alpha_p(t,S)
\end{multline}
where $D_A$ is the $|A|$-part of the discriminant $D$ (i.e. $D_A$
divides some power of $|A|$) and $H(m_1, Dt_2^2)=L(1-m_1, Dt_2^2)$ are
the H.~Cohen numbers. 
\end{theorem}

We should like to note that the variant of the Siegel formula given in
Theorem~\ref{rtS} is different from the formula given in~\cite{Sh}. 
Shimura used the $L$-function with a primitive character.  
We modify the local factors $\alpha_p(S)$ for the prime divisors of $|A|$
and use the function $L(s,\Delta)$ with a non-fundamental discriminant,
i.e., we put some other non-regular $p$-factors inside the $L$-function. 
As a result our
formulae (see Examples~\ref{ex:fivesquares}--\ref{ex:A5} below) are
shorter. 

\begin{proof}
From the definition of the local densities we see that
\begin{alignat}{3}
\alpha_p(t,S)=\alpha_p(2t,A)\quad&\text{ if }p\ne 2,\qquad&
\alpha_2(t,S)=2\alpha_2(2t,A). 
\end{alignat}

We assume that $p$ is not a divisor of $|A|$.  Let $l_p=\ord_p(t)$ and
$t=p^{l_p}t_{\bar p}$.  According to \cite[Hilfssatz~16]{Si} the
density $\alpha_p(2t, A)$ is given by 
\begin{align*}
\alpha_p(2t, A)=&
(1-p^{1-m})
\bigl(1+p^{2-m}+\dots+p^{(2-m)\frac{l_p-1}2}\bigr)\\
\intertext{for $l_p\equiv 1$ mod~$2$ and}
\alpha_p(2t, A)=&
(1-p^{1-m})
\biggl(1+p^{2-m}+\dots+p^{(2-m)(\frac{l_p}2-1)}+
\dfrac{p^{(2-m)\frac{l_p}2}}
{1-\varepsilon_{A,t}(p)p^{\frac{1-m}2}}\,\biggr)
\end{align*}
for $l_p\equiv 0$ mod~$2$,
where $\varepsilon_{A,t}(p)=\biggl(\dfrac{(-1)^{\frac{m-1}2}
|A|2t_{\bar p}}p\biggr)$.  If $l_p=0$ we take only the last summand in
the second bracket (see \cite[Hilfssatz 12]{Si}).  The numbers
$Dt_2^2$ and $2t|A|$ differ by a square $f^2$ such that $f$ divides
some power of $|A|$.  Therefore if $p$ does not divide $|A|$ and
$l_p$ is even then
\[
\varepsilon_{A,t}(p)=\chi_D(p)=\left(\sfrac{D}p\right)\ne 0. 
\] 
If $l_p$ is odd then $p|D$ and $\chi_D(p)=0$. 

Let us reorganise the $p$-factors in the Siegel formulae for the local
densities. We put
\begin{multline}
\alpha_p(t,S)=\left(\frac{1-p^{1-m}}{1-\chi_D(p)p^{\frac{1-m}2}}\right)
p^{(2-m)\lfloor\frac{l_p}2\rfloor}\\ 
\Bigl(1+\sum_{1\le j\le l_p/2}p^{(m-2)j}\Bigr)
\cdot\left(1-\chi_D(p)p^{\frac{1-m}2})\right). 
\end{multline}
This formula is valid for both even and odd $l_p$. If $l_p=1$
(i.e. if $p$ divides only $t_1$, not $t_2$) then $\alpha_p(t,S)=1-p^{1-m}$. 

Taking the product over all divisors of $t_2$ we obtain the factor
\[
t_2^{2-m}\sum_{d|t_2}d^{m-2}\prod_{p|d}(1-\chi_D(p)p^{\frac{1-m}2})=
\sum_{a|t_2}\mu(a)\chi_D(a)a^{\frac{1-m}2}\sigma_{2-m}\bigl(\sfrac{t_2}a\bigr). 
\]

Using the functional equation we can express $L(m_1, Dt_2^2)$ in terms
of $L(1-m_1, Dt_2^2)=H(m_1, Dt_2^2)$.  Together with the formula for
the Bernoulli numbers
\[
(-1)^{m_1+1}\frac{B_{2m_1}}{2m_1}
=\pi^{-\frac{1}2-2m_1}\Gamma(m_1)\Gamma(m_1+\sfrac{1}2)\zeta(2m_1)
\]
it gives us the second formula~\eqref{HCohen}.  We note that
$(-1)^{[m_1/2]}H(m_1,Dt_2^2)$ are positive rational numbers with
bounded denominators. The denominators are $120$ for $m_1=2$, 
$252$ for $m_1=3$, $240$ for $m_1=4$, etc.  (see \cite{C}). 
\end{proof}

The exact formulae for the local densities $\alpha_p(t,S)$,
$S(X)=\frac{1}2A[X]$, for all prime divisors of the determinant of $A$
including $p=2$ were calculated in many papers.  See for example
Malyshev~\cite{Mal}, who used a classical method of Gauss sums, and
Yang~\cite{Ya}, who calculated the local densities in terms of local
Whittaker integrals.  In the examples below we use the formulae of~\cite{Ya}. 

\begin{example}\label{ex:fivesquares}
The sum of five squares. 
\end{example} 
Let $S_5(X)=x_1^2+\dots+x_5^2$.  In this example we are finishing the
calculation of Siegel (see \cite[\S 10]{Si}) who found $r(t,S_{5})$
for odd $t$.  According to Theorem~\ref{rtS} we have
\[
r(t,S_{5})=t^{3/2}\frac{120}{\pi^2}L(2,Dt_2^2)
\frac{1-\chi_D(2)2^{-2}}{1-2^{-4}}\alpha_2(t,S_5),
\] 
where $t=2^at_1t_2^2$ with $a=2b$ or $2b+1$ as in Theorem \ref{rtS},
$D=\disc \QQ(\sqrt{t})$.  The formula for $\alpha_{2}(t,S)$ (see
\cite[pp. 323--324]{Ya}) is rather too long to give here. 
After some tedious transformations we obtain that
\[
\alpha_2(t,S_5)=1-\sum_{k=1}^b 2^{-3k+1}+(-1)^D2^{-3b-2}-\chi_D(2)2^{-3b-3},
\]
where $\chi_D(2)=0$ if $D\equiv 0\mod 4$ and $\chi_{D}(2)=1$ or $-1$
if $D\equiv 1 \mod 8$ or $D\equiv 5 \mod 8$ respectively.  In terms of
the Cohen numbers we have the numerical formula
\[
r(t,S_5)=
\begin{cases}
-40H(2,Dt_2^2)\cdot\dfrac{2^{3b+2}+3}7 &\text{if } D\equiv 0\mod 4,\\
\dfrac{-120H(2,Dt_2^2)}{4+\chi_D(2)}\cdot
\left(\dfrac{5\cdot 2^{3b+3}+2}7-\chi_D(2)\right)
&\text{if }D\equiv 1\mod 4. 
\end{cases}
\]
We note that $Dt_2^2$ is equal to $t$ up to a power of $2$. 

Let $L$ be an even integral quadratic lattice and $q_L(x)$ and $b_L(x,y)$
be the corresponding finite quadratic and bilinear forms on the
discriminant group $D(L)=L^\vee/L$. For $q_L$ we have the local
decomposition
\[
q_L=\bigoplus_p (q_L)_p=\bigoplus_p q_{L\otimes \ZZ_p},
\] 
where $q_{L\otimes \ZZ_p}$ is the finite quadratic form with values in
$\QQ_p/\ZZ_p$ ($p\ne 2$) or in $\QQ_2/2\ZZ_2$ and $(q_L)_p$ is the
discriminant form on the $p$-component of the finite abelian group
$L^\vee/L$ with values in $\QQ^{(p)}/\ZZ\cong \QQ_p/\ZZ_p$ or in
$\QQ^{(2)}/2\ZZ$ ($\QQ^{(p)}$ is the ring of fractions whose
denominator is a power of $p$).  A similar decomposition is valid for
$b_L$.  We recall that any quadratic form over the $p$-adic integers
$\ZZ_p$ ($p\ne 2$) is equivalent to a diagonal form. For $p=2$ it can
be represented as a sum of forms of types $2^nux^2$ ($u\in
\ZZ_2^*/(\ZZ_2^*)^2$), $2^n(2x_1x_2)$ and $2^n(x_1^2+2x_1x_2+x_2^2)$. 

\begin{example}\label{ex:A1plusD4}
The root lattice $A_1\oplus D_4$. 
\end{example} 
The quadratic form $S_{1,4}$ of this lattice is similar to the sum of
five squares. More exactly
\[
S_{1,4}=x_1^2+\frac{1}2(x_2^2+x_3^2+x_4^2+x_5^2),\quad
\text{where $x_2+x_3+x_4+x_5$ is even.} 
\] 
The determinant of $A_1\oplus D_4$ is equal to $8$.  The discriminant
form of $D_4$ is equal to the discriminant form of
$V(2)=2(2x_1^2+2x_1x_2+2x_2^2)$. Using this we obtain that over $\ZZ_2$
\[
\frac{1}2 (A_1\oplus D_4)\otimes{\ZZ_2}\cong\latt{2}\oplus 
\frac{1}2
\begin{pmatrix} 2&1\\1&2\end{pmatrix}\oplus
\begin{pmatrix} 2&1\\1&2\end{pmatrix}
\]
We use the notation of the previous example for $t=2^at_1t_2^2$,
$a=2b$ or $2b+1$.  Again using \cite{Ya} we obtain
\begin{equation}\label{a2AD}
\alpha_2(t,S_{1,4})=
1-\sum_{k=1}^b 2^{-3k}+(-1)^D2^{-3(b+1)}-\chi_D(2)2^{-3b-4}. 
\end{equation}
The second formula~\eqref{HCohen} of Theorem \ref{rtS} is 
\[
r(t,S_{1,4})=
\begin{cases}
-8H(2,Dt_2^2)\cdot 2^{3b+3}\alpha_2(t,S_{1,4}),& \quad\text{if }
D\equiv 0\mod 4,\\
\dfrac{-120H(2,Dt_2^2)}{4+\chi_D(2)}\cdot 2^{3b+4}\alpha_2(t,S_{1,4}),&
\quad \text{if }D\equiv 1\mod 4. 
\end{cases}
\]
The first formula~\eqref{L-zeta} of Theorem~\ref{rtS} gives us an
expression which we shall use later in our estimations of $N_{A_1\oplus
D_4}(2d)$:
\begin{equation}
r(t,S_{1,4})=t^{3/2}16\frac{\zeta(2)}{\zeta(4)}\,L(2,\, Dt_2^2)\,
\frac{1-\chi_D(2)2^{-2}}{1-2^{-4}}\,\alpha_2(t, S_{1,4}). 
\end{equation}
See \eqref{LsD} in order to understand the form of the factors. 

\begin{example}\label{ex:A5}
The root lattice $A_5$. 
\end{example} 
Let $D=\disc \QQ(\sqrt{3t})$ and $t=2^a3^ct_1t_2^2$. 
According to Theorem \ref{rtS}
\[
r(t,\sfrac{1}2A_5)=
t^{3/2}
\frac{32}{\sqrt{3}}\frac{\zeta(2)L(2,Dt_2^2)}{\zeta(4)}
\prod_{p=2,3}\frac{1-\chi_D(p)p^{-2}}{1-p^{-4}}
\alpha_p(t,\sfrac{1}2A_5). 
\]
The discriminant form of the lattice $A_5$ is the cyclic group of
order $6$ generated by the element $\bar v$ such that $\bar v\cdot
\bar v\equiv \frac{5}6\mod 2\ZZ$. For the local part of the
discriminant group we have
\[
\begin{array}{ll}
D(A_5)_3=\latt{2\bar v},\qquad &(2\bar v)^2\equiv \sfrac{1}3\mod \ZZ_3,\\
D(A_5)_2=\latt{3\bar v},\qquad  &(3\bar v)^2\equiv \sfrac{3}2\mod 2\ZZ_2. 
\end{array}
\]
It follows that 
\[
A_5\otimes \ZZ_3\cong x_1^2+x_2^2+x_3^2+2x_4^2+3x_5^2
\]
and
\[
A_5\otimes \ZZ_2\cong 2x_1x_2+2x_3x_4+6x_5^2\cong 2U\oplus \latt{6}. 
\]
Put $t'=2^at_1t_2^2$, so that $t=3^ct'$ and $(3,t')=1$.  The formula
for $\alpha_3$ (see \cite[p. 317]{Ya}) after some transformations can
be written as follows
\begin{equation}\label{a3A}
\alpha_3(3^ct',\sfrac 1{2}A_5)=
1-\sum_{k=1}^{3\lfloor\sfrac{c}2\rfloor+2}\biggl(\frac{k}3\biggr)3^{-k}
+\biggl(\frac{t'}3\biggr)3^{-\frac{3c+3}2},
\end{equation}
where $\bigl(\frac{k}3\bigr)$ is the Legendre symbol and we add the
last term only if $c$ is odd.  For $p=2$, we put $a=2b$ or
$2b+1$ and obtain
\begin{equation}\label{a2A}
\alpha_2(t, \sfrac{1}2A_5)=
1+\sum_{k=1}^b 2^{-3k-1}-(-1)^D2^{-3b-4}+\chi_D(2)2^{-3b-5}. 
\end{equation}
In terms of the Cohen numbers we have 
\[
r(2t, A_5)=
\left(\frac{t_{A_5}}{D_{A_5}}\right)^{3/2}\frac{1}{\sqrt{3}} 2^5\cdot 30(-H(2, Dt_2^2))
\prod_{p=2,3}\frac{1-\chi_D(p)p^{-2}}{1-p^{-4}}
\alpha_p(t,\sfrac{1}2 A_5),
\]
where $t_{A_5}=2^a3^c$ and $D_{A_5}$ are the products of the powers of $2$ and $3$ in
$t$ and $D$. 

\begin{proposition}\label{d56}
The inequality 
\[
30N_{A_1\oplus D_4}(2m)+16N_{A_5}(2m)<5N_{D_6}(2m)
\]
is true for any $m\ge 20$ and for $m=17$.  The inequality
\[
30N_{A_1\oplus D_4}(2m)+16N_{A_5}(2m)< 6N_{D_6}(2m). 
\]
is true if $m\ge 12$. 
\end{proposition}
\begin{proof} First we estimate $N_{D_6}(2m)$ from below.  By
  definition
\[
D_6=\{(x_i)\in \ZZ^6\mid x_1+\dots+x_6\in 2\ZZ\}. 
\]
Therefore the number $N_{D_6}(2m)$ is equal to the number of
representation of $2m$ by six squares. It is classically known (see
\cite[p. 187]{Iw}) and it can be easily proved using Eisenstein
series or the Siegel main formula that
\[
N_{D_6}(2m)=64\tilde\sigma_2(m,\chi_4)-4\sigma_2(m,\chi_4)
\]
where $\chi_4(m)=\bigl(\frac{-4}m\bigr)$ is the unique non-trivial
Dirichlet character modulo~$4$ and for any Dirichlet character $\chi$ we
put
\begin{alignat*}{2}
\sigma_k(m,\chi)=\sum_{d|m} \chi(d)d^k,
\qquad &
\tilde\sigma_k(m,\chi)=\sum_{d|m} \chi\left(\frac{m}d\right)d^k. 
\end{alignat*}
Let $a_p=\ord_p(m)$. For any quadratic character $\chi$ modulo~$\Delta$
we have
\begin{align*}
\tilde\sigma_k(m,\chi)&=m^k\sum_{p|m}
\frac{1-(\chi(p)p^{-k})^{(a_p+1)}}{1-\chi(p)p^{-k}}\\
&\ge m^k\prod_{p|m,\ (p,\Delta)=1}(1-p^{-k}). 
\end{align*}
This is because
\[
\tilde\sigma_k(m,\chi)=\sum_{d|m} \chi(d)\bigl(\frac{m}d\bigr)^k=
m^k\sum_{d|m} \chi(d)d^{-k}
\]
and
\[
\frac{1-(\chi(p)p^{-k})^{(a_p+1)}}{1-\chi(p)p^{-k}}
\ge 
\frac{1-p^{-k(a_p+1)}}{1+p^{-k}}
\ge 
\frac{1-p^{-2k}}{1+p^{-k}}=1-p^{-k}. 
\]

If $(m,\Delta)=1$ then
$\tilde\sigma_k(m,\chi)=\chi(m)\sigma_k(m,\chi)$ since
$\chi$ is a real character. 
Moreover for any prime divisor $p$ of the module
$\Delta$ of $\chi$
\begin{alignat*}{2}
\tilde\sigma_k(p^am_1,\chi)=p^{ak}\tilde\sigma_k(m_1,\chi),
\qquad & 
\sigma_k(p^am_1,\chi)=\sigma_k(m_1,\chi). 
\end{alignat*}
Therefore
\begin{equation}\label{ND6}
N_{D_6}(2m)\ge 60\tilde\sigma_2(m,\chi_{4})
> 60\zeta(2)^{-1}(1-2^{-2})^{-1} m^2= \frac{480}{\pi^2}\,m^2. 
\end{equation}

Next we have to estimate from above the Dirichlet series 
\[
\sum_{n\ge 1}\frac{b_n(\Delta)}{n^s}=
\frac{\zeta(s)L(s,\Delta)}{\zeta(2s)}
\] 
(see (\ref{LsD})) for $s=2$.  If $(n, \Delta)=1$ then
\[
b_n(\Delta)\le b_n(1)=2^{\rho(n)},
\]
where $\rho(n)$ is the number of prime divisors of $n$, with
equality if and only if $\bigl(\frac{\Delta}p\bigr)=1$ for any odd
prime divisor of $n$ and $\bigl(\frac{\Delta}8\bigr)=1$ if $n$ is
even.  If in $n$ there is at least one non-residue modulo $\Delta$
then $b_n(\Delta)=0$.  Therefore if $(p,\Delta)=1$ (it might be that
$p=2$) then the local $p$-factor of the Dirichlet series is equal to
\begin{equation}\label{goodp}
1+2\sum_{m\ge 1}p^{-ms}=\frac{p^s+1}{p^s-1}. 
\end{equation}
Let us assume that $\Delta=p^{2k}\Delta'$ (it might be that $p=2$)
with $(p,\Delta')=1$.  Considering the congruence class of
$b_{p^m}(p^{2k}\Delta')$ for all powers of $p$ we see that the local
$p$-factor of the Dirichlet series equals
\begin{equation}\label{badp}
1+\sum_{m=1}^{2k}\frac{p^{\lfloor\frac{m}2\rfloor}}{p^{ms}}+
2p^k\sum_{m\ge 2k+1}p^{-ms}. 
\end{equation}
If $\Delta=p^{2k+1}\Delta'$ then the local factor is smaller: the last
term in~\ref{badp}, $2p^k\sum_{m\ge 2k+1}p^{-ms}$, is replaced by one
summand, $p^{k-(2k+1)s}$. A direct calculation shows that for $s=2$ the
regular factor (\ref{goodp}) is larger than the non-regular factor
(\ref{badp}) for any prime $p\ge 2$.  Therefore
\[
\frac{\zeta(2)L(2,\Delta)}{\zeta(4)}\le \prod_{p}
\frac{p^2+1}{p^2-1}=\frac{\zeta(2)^2}{\zeta(4)}=\frac{5}2. 
\]
The next step is an estimation from above of the $2$-factor in
$N_{A_1\oplus D_4}(2m)=r(m,S_{1,4})$ and the $2$- and $3$-factors in
$N_{A_5}(2m)=r(m,\frac{1}2A_5)$. Elementary calculation using
(\ref{a2AD}) gives us
\begin{equation}\label{2factor_A+D}
\frac{1-\chi_D(2)2^{-2}}{1-2^{-4}}\,\alpha_2(m, S_{1,4})\le \frac{5}4, 
\end{equation}
with equality if $D\equiv 5\mod 8$ and $m$ is odd. 

For the local $3$-factor in $\sfrac{1}{2}A_5$ we obtain (using~(\ref{a3A}))
\begin{equation}\label{3-factor_A5}
\frac{1-\chi_D(3)3^{-2}}{1-3^{-4}}\,\alpha_3(m, \sfrac{1}2A_5)\le \frac{11}{12}  
\end{equation}
with equality if $m=3m'$, where $m'\equiv 2\mod 3$. 

For the local $2$-factor in $\sfrac{1}{2}A_5$ we obtain
\begin{equation}\label{2factor_A5}
\frac{1-\chi_D(2)2^{-2}}{1-2^{-4}}\,\alpha_2(m,
\frac{1}2A_5)\le \frac{10}{7}.
\end{equation}
In this case we must analyse the case when $m=2^am'$ and
$b=\lfloor\sfrac{a}2\rfloor$ goes to infinity (see (\ref{a2A})).  If
$D\equiv 0\mod 4$ or $\equiv 5\mod 8$ then the local density tends to
its supremum as $b$ tends to infinity.  This value is equal to
$\frac{15}{14}$.  Therefore the left-hand side of~(\ref{2factor_A5})
is smaller than $\frac{10}7$.  If $D\equiv 1\mod 8$ then $\alpha_2$
takes its maximal values $\frac{35}{32}$ for $b=0$. In this case
the left-hand side of~(\ref{2factor_A5}) is equal to $\frac{7}8$. 

Now we can combine all our estimates.  We have
\begin{alignat*}{2}
N_{A_1\oplus D_4}(2m)\le 50 m^{3/2},
\qquad &
N_{A_5}(2m)\le \frac{2200}{21\sqrt{3}} m^{3/2}. 
\end{alignat*}
Using (\ref{ND6}) we obtain that the inequalities~(\ref{mineq1}) and
(\ref{mineq2}) of Theorem~\ref{mainineq} and Proposition~\ref{d56}
are valid for $m> 102$ and for $m>71$ respectively.  For
many $m$ smaller or equal to $102$ we can write a better estimate
for the number of representations.  But we can use the exact formulae
for the theta series of $D_6$, $A_1\oplus D_4$ and $A_5$ in terms of
Jacobi theta series in order to check the inequality for $m\le 102$. 

The theta series of the lattice $A_n$ is given by the formula
(see \cite[Ch. 4, (56)]{CS})
\[
\theta_{A_5}(\tau)=
\dfrac{\sum_{k=0}^5 \vartheta_3(\tau, \frac{k}{6})^{6}}
{6\vartheta_3(6\tau)},
\]
where
\[
\vartheta_3(\tau,z)=
\sum_{n\in \ZZ}\exp(\pi i(n^2\tau+2nz)),
\]
and
$\vartheta_3(\tau)=\vartheta_3(\tau,0)$. 
For the lattice $D_n$ one has (see \cite[Ch. 4, (87), (10)]{CS}) that
\[
\theta_{D_n}(\tau)
=\frac{1}{2}(\vartheta_3(\tau)^n+\vartheta_3(\tau+1)^n). 
\]
Using these formulae we can compute the first $102$
Fourier coefficients of the function
\[
5\theta_{D_6}-30\theta_{A_1\oplus D_4}-16\theta_{A_5}. 
\]
(This is a calculation with formal power series which we did using PARI.)  
We find that these coefficients are negative exactly for $d < 20$ and $d \neq
17$. Hence the first inequality of the proposition holds as stated. Repeating
the same calculation with $6$ instead of $5$ we obtain the second inequality. 
\end{proof}

We have now proved a slightly weakened version of
Theorem~\ref{gentype}. To obtain the full result we need the following
observation. 

\begin{proposition}\label{search}
For $d=12$, $13$, $14$, $15$, $16$, $18$ and $19$ there exist vectors
$l_d\in E_7$ that satisfy $l_d^2=2d$ that are orthogonal to at least
$2$ and at most $14$ roots. For $d=9$ and $d=11$ there exist vectors
of length $l^2=2d$ that are orthogonal to exactly $16$ roots. 
\end{proposition}
\begin{proof}
These were found by a computer search. We give one example in each
case. We express the vectors in terms of the simple roots $v_i$, $1\le
i\le 7$ which are given in terms of the standard basis
$e_1,\ldots,e_8$ of $\QQ^8$ by
\begin{align*}
v_i&=e_{i+2}-e_{i+1}\text{\quad for }1\le i\le 6,\\
v_7&=\sfrac{1}{2}(e_1+e_2+e_3+e_4)-\sfrac{1}{2}(e_5+e_6+e_7+e_8)
\end{align*}
(see~\cite{Bou}). The examples are shown in
Table~\ref{ta:shortvectors}: the vector $l_d=\sum \lambda_{d,i} v_i$
with $\lambda_d=(\lambda_{d,1},\ldots,\lambda_{d,7})\in\ZZ^7$ is
orthogonal to exactly $2p_d$ roots of $E_7=\sum_{i=1}^7\ZZ
v_i\subset\QQ^8$. 
\begin{table}
\begin{center}
\begin{tabular}{|c|c|c|}
\hline
$d$&$p$&$\lambda_d$\\
\hline
 9 & 8 & -1,2,3,1,2,1,3\\ \hline
11 & 8 & 3,3,0,-1,-2,-1,0\\ \hline
12 & 7 & 2,1,2,-2,0,0,1\\ \hline
13 & 7 & 2,3,-1,1,0,0,-1\\ \hline
14 & 6 & 2,0,3,0,2,1,1\\ \hline
15 & 7 & 1,-2,0,2,4,2,0\\ \hline
16 & 6 & 1,0,-1,3,0,0,-2\\ \hline
18 & 5 & 3,2,3,2,0,0,-2\\ \hline
19 & 6 & 2,3,2,-3,-4,-2,1\\ \hline
\end{tabular}
\caption{Short vectors in $E_7$ orthogonal to few roots}\label{ta:shortvectors}
\end{center}
\end{table}
There are other vectors with the required properties (for instance, we
found one with $d=19$ and $p=7$), but none for smaller~$d$. 
\end{proof}

\bibliographystyle{alpha}

\bigskip
\noindent
V.A.~Gritsenko\\
Universit\'e Lille 1\\
Laboratoire Paul Painlev\'e\\
F-59655 Villeneuve d'Ascq, Cedex\\
France\\
{\tt valery.gritsenko@math.univ-lille1.fr}
\bigskip

\noindent
K.~Hulek\\
Institut f\"ur Algebraische Geometrie\\
Leibniz Universit\"at Hannover\\
D-30060 Hannover\\ 
Germany\\
{\tt hulek@math.uni-hannover.de}
\bigskip

\noindent
G.K.~Sankaran\\
Department of Mathematical Sciences\\
University of Bath\\
Bath BA2 7AY\\
England\\
{\tt gks@maths.bath.ac.uk}

\end{document}